\documentclass[12pt,a4paper,twoside]{amsart}
\usepackage[UKenglish]{babel}
\usepackage[T1]{fontenc}
\usepackage[utf8x]{inputenc}
\usepackage{latexsym,amsfonts,amsmath,amsthm,amssymb,mathrsfs}

\usepackage{scrtime}
\usepackage[colorlinks]{hyperref}

\usepackage{fullpage}
\usepackage{xcolor}
\usepackage{paralist}
\usepackage{todonotes}
\usepackage{dsfont}
\usepackage{float}
\usepackage{caption}
\usepackage{xfrac}
\usepackage{algorithm,algpseudocode}
\usepackage[pagewise]{lineno}

\usepackage{graphicx}

\usepackage{enumitem}

\newlist{enumth}{enumerate}{1}
\setlist[enumth]{label=\emph{(\arabic*)}, ref=(\arabic*)}

\usepackage{fourier}

\newcommand{\aaa}{c}
\newcommand{\cS}{{\mathcal S}}
\renewcommand{\subseteq}{\subset}

\DeclareMathOperator{\R}{\mathbf{R}}
\DeclareMathOperator{\N}{\mathbf{N}}

\DeclareMathOperator{\nextp}{\mathrm{next}}

\newtheorem*{rem*}{Remark}

\DeclareMathSymbol{A}{\mathalpha}{operators}{`A}%
\DeclareMathSymbol{B}{\mathalpha}{operators}{`B}%
\DeclareMathSymbol{C}{\mathalpha}{operators}{`C}%
\DeclareMathSymbol{D}{\mathalpha}{operators}{`D}%
\DeclareMathSymbol{E}{\mathalpha}{operators}{`E}%
\DeclareMathSymbol{F}{\mathalpha}{operators}{`F}%
\DeclareMathSymbol{G}{\mathalpha}{operators}{`G}%
\DeclareMathSymbol{H}{\mathalpha}{operators}{`H}%
\DeclareMathSymbol{I}{\mathalpha}{operators}{`I}%
\DeclareMathSymbol{J}{\mathalpha}{operators}{`J}%
\DeclareMathSymbol{K}{\mathalpha}{operators}{`K}%
\DeclareMathSymbol{L}{\mathalpha}{operators}{`L}%
\DeclareMathSymbol{M}{\mathalpha}{operators}{`M}%
\DeclareMathSymbol{N}{\mathalpha}{operators}{`N}%
\DeclareMathSymbol{O}{\mathalpha}{operators}{`O}%
\DeclareMathSymbol{P}{\mathalpha}{operators}{`P}%
\DeclareMathSymbol{Q}{\mathalpha}{operators}{`Q}%
\DeclareMathSymbol{R}{\mathalpha}{operators}{`R}%
\DeclareMathSymbol{S}{\mathalpha}{operators}{`S}%
\DeclareMathSymbol{T}{\mathalpha}{operators}{`T}%
\DeclareMathSymbol{U}{\mathalpha}{operators}{`U}%
\DeclareMathSymbol{V}{\mathalpha}{operators}{`V}%
\DeclareMathSymbol{W}{\mathalpha}{operators}{`W}%
\DeclareMathSymbol{X}{\mathalpha}{operators}{`X}%
\DeclareMathSymbol{Y}{\mathalpha}{operators}{`Y}%
\DeclareMathSymbol{Z}{\mathalpha}{operators}{`Z}%

\renewcommand{\leq}{\leqslant}
\renewcommand{\geq}{\geqslant}

\renewcommand{\ge}{\geqslant}
 
\numberwithin{equation}{section}

\newcommand{\uple}[1]{\text{\boldmath${#1}$}}
\newcommand{\bm}{\uple}

\def\setminus{\mathchoice
    {\mathbin{\vrule height .62ex width 1.61ex depth -.38ex}}
    {\mathbin{\vrule height .62ex width 1.61ex depth -.38ex}}
    {\mathbin{\vrule height .50ex width 0.85ex depth -.28ex}}
    {\mathbin{\vrule height .20ex width 0.570ex depth -.24ex}}
}

\renewcommand{\mathcal}{\mathscr}

\newcommand{\fl}[1]{\lfloor {{#1}}\rfloor}

\newcommand{\Cc}{\mathbf{C}}

\newcommand{\Zz}{\mathbf{Z}}

\newcommand{\Rr}{\mathbf{R}}

\newcommand{\Qq}{\mathbf{Q}}

\newcommand{\expect}{\mathbf{E}}

\def\loccit{loc.\kern3pt cit.{}\xspace}
\def\cf{see\kern.3em}
\def\Cf{See\kern.3em}
\def\eg{e.g.\kern.3em}

\def\resp{\text{resp.}\kern.3em}

\newcommand{\mods}[1]{\,(\mathrm{mod}\,{#1})}

\newcommand{\eps}{\varepsilon}
\renewcommand{\rho}{\varrho}

\newcommand{\demi}{{\textstyle{\frac{1}{2}}}}

\theoremstyle{plain}
\newtheorem{theorem}{Theorem}[section]
\newtheorem*{theorem*}{Theorem}
\newtheorem{lemma}[theorem]{Lemma}
\newtheorem{corollary}[theorem]{Corollary}

\newtheorem{proposition}[theorem]{Proposition}

\theoremstyle{remark}

\theoremstyle{definition}

\newtheorem{definition}[theorem]{Definition}

\newtheorem{remark}[theorem]{Remark}

\allowdisplaybreaks
\title{Rational approximation with chosen numerators}
\begin{document}
\date{}

\author{Manuel Hauke}
\address[M. Hauke]{Dep. of Math. Sciences, NTNU Trondheim,
Sentralbygg 2, 7034 Trondheim, Norway\\\phantom{(M. Hauke)}Inst. of Analysis and Number Theory, TU Graz, Steyrergasse 30/II, 8010 Graz, Austria}
\email{hauke@math.tugraz.at}

\author{Emmanuel Kowalski}
\address[E. Kowalski]{D-MATH, ETH Z\"urich, R\"amistrasse 101, 8092 Z\"urich, Switzerland} 
\email{kowalski@math.ethz.ch}

\begin{abstract}
  We consider the problem of approaching real numbers with rational
  numbers with prime denominator and with a single numerator allowed
  for each denominator. We obtain basic results, both probabilistic
  and deterministic, draw connections to twisted diophantine approximation, and present a simple application, related to
  possible correlations between trace functions and dynamical
  sequences.
\end{abstract}

\makeatletter
\@namedef{subjclassname@2020}{%
  \textup{2020} Mathematics Subject Classification}
\makeatother


\maketitle

\section{Introduction}

Classical diophantine approximation, starting from Dirichlet's
Theorem, aims at approaching (irrational) real numbers $x\in [0,1]$
with rational numbers $a/q$, where~$q\geq 1$ is an integer
and~$(a,q)=1$, so that the distance $|x-a/q|$ is as small as possible
in terms of the height~$q$.

Motivated by a question at the interface of ergodic theory and the
theory of trace functions (as explained in
Section~\ref{sec-application}), we consider in this paper a different
approximation problem: first, we only consider prime denominators~$p$,
and more importantly, for each prime number~$p$, we assume that a
\emph{single} numerator~$a_p$, $0\leq a_p<p$, has been chosen. We
further fix a positive real number~$c>0$ and ask, for a given~$x$ or
for almost all~$x$, whether the inequality
$$
\Bigl|x-\frac{a_p}{p}\Bigr|\leq \frac{c}{p}
$$
has infinitely many solutions or not.

Clearly, the answer to such questions depends on the choice of the
sequence $(a_p)$ and is sometimes (trivially) negative, but we will see that it is also
sometimes positive in various cases. To state the results precisely,
we introduce the following notation. We let
$$
\Omega=\{(a_p)_p \,\mid\, a_p\in\Zz,\ 0\leq a_p<p\text{ for all } p\}
$$
(where, here and below, an index~$p$ indicates a family or sum over
primes). This is a compact Hausdorff space when each finite set
$\{0,\ldots,p-1\}$ is given the discrete topology. For
$\bm{a}\in\Omega$ and~$c>0$, we define
\[
  \mathcal{A}(\bm{a},c) = \Bigl\{ x \in [0,1]\,\mid\, \Bigl|
  x - \frac{a_p}{p}\Bigr|\leq \frac{c}{p}\text{ for infinitely
    many } p\Bigr\},
\]
and
\[
  \mathcal{A}(\bm{a}) = \bigcap_{c>0}\mathcal{A}(\bm{a},c) = \Bigl\{ x
  \in [0,1]\,\mid\, \liminf_{p \to \infty} p\Bigl| x -
  \frac{a_p}{p}\Bigr| = 0\Bigr\}.
\]

We denote by~$\lambda$ the Lebesgue measure on~$[0,1]$ and further
endow~$\Omega$ with the probability measure where each
coordinate~$a_p$ is uniformly distributed in $\{0,\ldots,p-1\}$. We
then have a first introductory result, which is probabilistic in nature:

\begin{theorem}\label{th-random}
  For almost all~$\bm{a}\in \Omega$, we have
  \[
    \lambda(\mathcal{A}(\bm{a}))=1.
  \]

  In fact, for fixed~$c\in ]0,\demi]$, and for almost
  all~$\bm{a}\in\Omega$, there exists a set~$S\subset [0,1]$ of
  Lebesgue measure~$1$ such that for~$x\in S$, for~$\eps>0$, and
  for~$X\geq 2$, we have
  \begin{equation}\label{eq-count}
    |\{p\leq X\,\mid\, |x-a_p/p|\leq c/p\}|=
    2c\log\log X+O((\log\log X)^{1/2+\eps}).
  \end{equation}
\end{theorem}

As a uniform version in $c$, we recover the following result.

\begin{corollary}\label{uniform_c}
  For almost every $\bm{a} \in\Omega$, there exists a set
  $S\subset [0,1]$ of full Lebesgue measure such that, for any
  $0 <c < 1/2$ and any $x \in S$, we have
  \[
    \Bigl|\Bigl\{p \leq X\,\mid\,
    \bigl\lvert x - \tfrac{a_p}{p} \bigr\rvert
    \leq \tfrac{c}{p}\Bigr\}\Bigr|
    \sim 2c\log \log X, 
  \]
  as $X \to \infty$.
\end{corollary}

While at least the almost sure divergence might not come as a surprise (see e.g. \cite{R20} for a related concept without the prime restriction), we will present different
proofs, in particular to point out an interesting analogy with sieve
which does not seem to be well-known.\\

Although Theorem~\ref{th-random} proves the existence of sequences
$\bm{a}\in\Omega$ with $\lambda(\mathcal{A}(\bm{a}))=1$, which is
sufficient for the application in Section~\ref{sec-application}, it is
natural to want to exhibit \emph{concrete} examples of sequences with
this property. We show that this can be done in a variety of ways.

The most obvious way of trying to find a suitable sequence is to apply
a greedy algorithm, trying at each step to cover the (essentially)
maximal space of the unit interval which has not been covered
recently. This turns out to be successful.

\begin{theorem}\label{thm_greedy}
  The deterministic greedy algorithm presented in detail in
  Section~\textup{\ref{sec-greedy}} generates a sequence
  $\bm{g}\in\Omega$ such that $\lambda(\mathcal{A}(\bm{g})) = 1$.  In
  fact, this sequence has the property that
  $\mathcal{A}(\bm{g},c)=[0,1]$ if~$c>2$, whereas for $c>0$ sufficiently
  small, we have
  \[
    \dim_H([0,1] \setminus \mathcal{A}(\bm{g},c)) = 1.
  \]
\end{theorem}

The sequence~$\bm{g}$ satisfies~$\lambda(\mathcal{A}(\bm{g}))=1$, but
it is very special, and, in fact has certain undesirable, or at least
unusual, properties. For instance, one can show that the sequence
$(g_p/p)$ is not uniformly distributed in $[0,1]$ (with respect to
Lebesgue measure), whereas almost all~$\bm{a}\in\Omega$ satisfy this
property.

Our next examples are more satisfactory in that respect. Describing
them requires some further notation. First, for a prime number~$p$, we
denote by $n_p\geq 1$ the positive integer such that $p$ is the
$n_p$-th prime in increasing order (in other words, we have
$n_p=\pi(p)$, where~$\pi(x)$ is the usual prime-counting
function). Moreover, if~$p$ is a prime and~$n\in\Zz$, we denote
by~$n\pmod p$ the unique integer congruent to~$n$ modulo~$p$ in the
set $\{0,\ldots,p-1\}$.

We will be considering sequences of the form $\uple{a}(\beta)$
and~$\uple{b}(\beta)$, with
\[
  a_p(\beta)=\lfloor \beta pn_p\rfloor \pmod p,\quad\quad
  b_p(\beta)=\lfloor \beta p^2\rfloor \pmod p
\]
for some real number~$\beta>0$. We will show that, under suitable Diophantine conditions on $\beta$, these satisfy the property
$\lambda(\mathcal{A}(\uple{a}))=\lambda(\mathcal{A}(\uple{b}))=1$.\\

We recall that a number $\beta \in \R$ is called \textit{badly approximable} when there exists $c = c(\beta) > 0$ such that
$\left\lvert \beta - \frac{a}{q}\right\rvert > \frac{c}{q^2}$ for all $a \in \mathbb{Z}, q \in \N$.
The set of badly approximable numbers contains all real algebraic numbers of degree~$2$ such as the Golden Ratio or $\sqrt{2}$, has Lebesgue measure~$0$, but Hausdorff dimension~$1$, and coincides
with the set of numbers whose partial quotients are bounded.\\

While our statement is more general, it contains the following simple result as an important special case:
\begin{theorem}\label{thm_badly}
  If~$\beta$ is badly approximable, then the sequences
  $\uple{a}(\beta)$ and~$\uple{b}(\beta)$ above satisfy
  \[
    \lambda(\mathcal{A}(\bm{a}(\beta)))=\lambda(\mathcal{A}(\bm{b}(\beta))
    = 1.
  \]
\end{theorem}

We define the Diophantine exponent of $\beta$ to be
\[\kappa = \kappa(\beta) := \sup \left\{\tau > 0: \left\lvert \beta - \frac{a}{q}\right\rvert < \frac{1}{q^{1 +\tau}} \text{ for i.m. } (a,q) \in \mathbb{Z} \times \N\right\}.\]
The set of irrationals with Diophantine exponent $1$ forms a set of full Lebesgue measure and by Roth's Theorem contains any irrational real algebraic number. However, for every $\kappa \geq 1$, there exist irrational numbers with Diophantine exponent $\kappa$, and even examples with $\kappa = \infty$ (the Liouville numbers) can be constructed.

As already mentioned, the statement we prove applies to more general
numbers~$\beta$, including some with arbitrary Diophantine exponent. To state the precise result requires the following 
diophantine condition.

\begin{definition}\label{iba_def}(Intermittently badly approximable numbers)
  Let~$\beta\in ]0,1[$ be a an irrational number with continued
  fraction expansion
  \[
    \beta = [0;a_1,a_2,\ldots].
  \]

  For~$k\geq 1$, let
  \[
    \frac{p_k}{q_k} = [0;a_1,a_2,\ldots,a_k].
  \]

  We say that~$\beta$ is \emph{intermittently badly approximable},
  abbreviated i.b.a, if there exist a real number $\delta$ with
  $0<\delta<1$, an integer~$B\geq 1$ and increasing sequences of
  integers $(j_k)_{k \geq 1}$ and $(T_k)_{k \geq 1}$ such that
  \begin{gather}
      \label{bad_range}
      a_{j_k}\leq B,\quad\ldots\quad,a_{j_k+T_k}\leq B,
      \\
    \label{range_length}
    T_k \geq \delta \log q_{j_k}
  \end{gather}
  for all~$k\geq 1$.
\end{definition}

\begin{remark}
  It is immediate that any badly approximable number~$\beta$ is
  intermittently badly approximable.
  However, since the definition allows very good approximations
  between the ``stretches of bad approximation'', one can easily check
  that there exist i.b.a. numbers of every
  possible (even infinite) Diophantine exponent.
  We remark that since the set of i.b.a. numbers contains all badly approximable numbers, the set has full Hausdorff dimension.
On the other hand, the set of i.b.a. numbers has Lebesgue measure $0$: To see this, we use a classical result of Diamond and Vaaler \cite{dv86} stating that for (Lebesgue) almost every $\beta = [0;a_1,a_2,\ldots]$, we have
\begin{equation}\label{dv_conv}\sum_{\ell \leq K}a_{\ell} - \max_{\ell \leq K}a_{\ell} \sim \frac{K \log K}{\log 2}.\end{equation}
Assuming an i.b.a. $\beta$ to satisfy \eqref{dv_conv}, we choose $K_1 = K_1(k) = j_k, K_2 = K_2(k) = j_k + T_k$ with $j_k,T_k,B$ as in \eqref{bad_range} and \eqref{range_length}.
Note that this would imply

\begin{equation}\label{impossible_eq}\begin{split}\frac{(K_1 + T_k)\log K_1}{2} &\leq \frac{K_2 \log K_2}{\log 2} \sim \sum_{\ell \leq K_2}a_{\ell} - \max_{\ell \leq K_2}a_{\ell}
\\&\leq T_k\cdot B + \sum_{\ell \leq K_1}a_{\ell}  - \max_{\ell \leq K_1}a_{\ell} \sim T_k\cdot B + \frac{K_1 \log K_1}{\log 2}.
\end{split}
\end{equation}

Since $(q_k)_{k \in \mathbb{N}}$ grows at least exponentially, \eqref{range_length} implies that $T_k \gg_{\delta} j_k = K_1$, so for $K_1$ sufficiently large, \eqref{impossible_eq} becomes impossible.
\end{remark}

Definition \ref{iba_def} allows us now to state Theorem \ref{thm_badly} in a more general fashion.

\begin{theorem}\label{th-deterministic}
  Let~$\beta$ be an i.b.a. number.
  Then the sequences $\bm{a}(\beta) = (a_p)_{p}$ and~$\bm{b}(\beta) = (b_p)_{p}$
  defined by
  \[
    a_p = \lfloor pn_p\beta\rfloor \pmod p,\quad\quad b_p = \lfloor
    p^2\beta \rfloor \pmod p
  \]
  satisfy $\lambda(\mathcal{A}(\bm{a}(\beta))) = \lambda(\mathcal{A}(\bm{b}(\beta))) = 1$.
\end{theorem}

In the converse direction, we will show that some Diophantine condition needs to be assumed:

\begin{theorem}
  There exist irrational real numbers~$\beta$ such that
  $\lambda(\mathcal{A}(\bm{a}(\beta))) = \lambda(\mathcal{A}(\bm{b}(\beta))) = 0$.
\end{theorem}

The proofs of these results are found in
Section~\ref{sec-deterministic}.\\

As in the Greedy sequence, we also obtain some information on the Hausdorff dimension of the
exceptional sets, showing that we cannot replace ``full measure'' with ``the whole unit interval'':

\begin{theorem}\label{th-hd}
  For~$\beta$ a badly approximable real number, we have
  \[
    \dim_H([0,1] \setminus \mathcal{A}(\bm{a}(\beta)) )= \dim_H([0,1]
    \setminus \mathcal{A}(\bm{b}(\beta)) ) = 1.
  \]
\end{theorem}

See Section~\ref{sec-HD} for the proof.

\begin{remark}
  (1) A natural question is whether the result extends to sequences
  generated by (say, monic) polynomials $f\in\Zz[X]$, i.e., when
  $a_p = \fl{\beta p f(n_p)}\mods{p}$ or $\fl{\beta p f(p)}\mods{p}$
  (at least when $\beta$ is badly approximable).

  (2) Another concrete sequence for which an answer would be interesting
  is the following, where we use a very slight variant of our
  setting. For $p\equiv 1\mods{4}$, let $a_p$ and $b_p$ be the
  representatives in $\{0,\ldots,p-1\}$ of the solutions of the equation
  $n^2\equiv -1\mods{p}$. Is it true (or not) that for almost
  all~$x\in [0,1]$, there are infinitely many primes $p\equiv 1\mods{4}$
  such that either $|x-a_p/p|<c/p$ or $|x-b_p/p|<c/p$?  This question is
  of course related to the celebrated result of Duke, Friedlander, and
  Iwaniec, according to which the roots of $X^2+1=0$ modulo primes
  $p\leq x$ are equidistributed in $[0,1]$ for the Lebesgue measure as
  $x\to+\infty$ (see~\cite[\S\,21.3]{ant} for a proof).

  (3) We expect that the sequences in Theorem~\ref{th-deterministic}
  should satisfy the second generic property of
  Theorem~\ref{th-random}, at least for $\beta$ badly approximable and as far as the main term is
  concerned. In other words, we expect that if $0<c<1/2$, then there
  exists
  a set $S \subseteq [0,1]$ of Lebesgue measure~$1$ such that, for any
  $\alpha \in S$, we have
  \begin{gather*}
    |\{p \leq Q\,\mid\, | \alpha - a_{p}(\beta)/p|\leq c/p\}|
    \sim 2c\log \log Q,
    \\
    |\{p \leq Q\,\mid\, | \alpha - b_{p}(\beta)/p|\leq c/p\}|
    \sim 2c\log \log Q,
  \end{gather*}
  as~$Q\to +\infty$.
\end{remark}

It turns out that the above statements connect quite smoothly (see Section \ref{sec-deterministic}) to the topic of \textit{twisted diophantine
	approximation}, a research area that originates in the work of Kurzweil \cite{K55}.
Conversely to classical (inhomogeneous) metric Diophantine approximation, where usually the inhomogeneous parameter is fixed while the rotation parameter is randomized, here the rotation parameter $\beta$ is fixed, and one asks about the measure of the set \[
A_{\beta}(\psi) = \left\{x \in [0,1)\,\mid\, \lVert n\beta - x \rVert \leq
\psi(n) \text{ for infinitely many integers } n\geq 1\right\}.
\]
where $\psi: \N \to [0,\infty[$ is some approximation function.\\
While we will use Kurzweil's Theorem to deduce the statement of Theorem \ref{th-deterministic} for $a_p = \lfloor pn_p\beta\rfloor \pmod p$, we will use the connection in the other direction to deduce from Theorem \ref{th-deterministic} for $b_p = \lfloor p^2\beta \rfloor \pmod p$ the following result in twisted, denominator-restricted Diophantine approximation, that might be of independent interest:

\begin{theorem}
	\label{prime_twisted}
	Let~$\beta$ be an irrational real number which is intermittently badly approximable. Then
	\[
	\lambda\Bigl(\bigl\{x \in [0,1)\,\mid\, \liminf_{p \to \infty} p\lVert
	p\beta - x \rVert = 0\bigr\}\Bigr)=1.
	\]
\end{theorem}

\begin{remark}
    We remark that Theorem \ref{prime_twisted} also holds for Lebesgue almost all $\beta$ instead of $\beta$ being i.b.a. This follows from the (straightforward) doubly metric inhomogeneous Khintchine Theorem due to Cassels \cite[Chapter VII, Theorem II]{Cas57}): For arbitrary $\psi: \mathbb{N}: \to [0,1/2]$ with $\sum_{n \in \mathbb{N}} \psi(n) = \infty$, we have (with $\lambda_2$ denoting the two-dimensional Lebesgue measure)
    \[\lambda_2\left(\left\{(\beta,x) \in [0,1)^2: \lVert n\beta - x\rVert \leq \psi(n)\text{\;i.o.}\right\}\right) = 1.\]
    Then choosing $\psi(n) = \frac{\varepsilon}{n}\mathds{1}_{[n \in \mathbb{P}]}$, applying Fubini's Theorem and $\varepsilon \to 0$ immediately proves the claim.
    However, note that in this result we do not have any information about the full measure set in $\beta$ at all, and therefore cannot provide any \textit{explicit} example of such $\beta$, in contrast to the result of Theorem \ref{prime_twisted}. 
\end{remark}

\subsection*{Notation}

We use the usual $O$ and $o$-notation as well as Vinagradov notation
$\gg$ and~$\ll$. If $f(x) \ll g(x)$ and~$g(x) \ll f(x)$ for all~$x$
in a set~$X$ on which~$f$ and~$g$ are defined, we write
$f(x) \asymp g(x)$ on~$X$. We write $f\sim g$ as $x\to x_0$ if
$\lim_{x\to x_0} f(x)/g(x)=1$.

We denote $e(x)=e^{2\pi ix}$.
We write the indicator function for a set $A$ as $\mathds{1}_A$.

Given a set $X$ and a sequence $(A_n)_{n\geq 1}$ of subsets of~$X$, we
denote
$$
\limsup_{n}A_n=\bigcap_{N\geq 1}\bigcup_{M\geq N}A_M,
$$
which is the set of elements of~$X$ that belong to infinitely many
sets~$A_n$.

\subsection*{Acknowledgements}

A part of this work was supported by
the Swedish Research Council under grant no. 2016-06596 while MH and EK were in residence at Institut
Mittag-Leffler in Djursholm, Sweden in 2024. Another essential part of this work was carried out while MH was visiting EK at ETH Z\"urich. MH would like to thank the university for the support and hospitality.
 MH would also like to thank Christoph Aistleitner, Victor Beresnevich, Sam Chow, Julia Stadlmann, and Sanju Velani for various discussions on related topics. Further, we thank A. Zafeiropoulos and N. Moshevitin for pointing out to us the
relevance of the theory of twisted diophantine approximation to the
problem considered here. Further, we would like to thank the anonymous reviewers for their valuable and pertinent comments that clearly improved the readability of the manuscript.

\section{Application to trace functions}
\label{sec-application}

This section is independent of the rest of the paper, but it explains
the original motivation for the study of the particular approximation
problem that we consider; more details and discussion can be found
in~\cite{kowalski_ergodic}. The question is whether certain types of
``triangular'' ergodic averages converge almost everywhere, and we show
that this is sometimes only the case when the limit is taken along
sufficiently sparse subsequences. 

Precisely, let $X=(\Rr/\Zz)^2$ and $\mu$ the (probability) Lebesgue measure
on~$X$. Further, let $f\colon X\to X$ be the map defined by
$f(x,y)=(x+y,y)$. We have $f_*\mu=\mu$. Define
$\varphi\colon X\to \Cc$ by $\varphi(x,y)=e(x)$.

We chose a sequence $(a_p)$ as in Theorem~\ref{th-random}, with
$c=1/2$ (e.g., one of the explicit sequences given by
Theorem~\ref{th-deterministic}, say
$a_p=\lfloor \sqrt{2}p^2\rfloor\mods{p}$). For $p$ prime and
$n\in\Zz$, we define $t_p(n)=e(-na_p/p)$.

\begin{proposition}
  For $p$ prime, define $s_p\colon X\to \Cc$ by
  $$
  s_p(x,y)=\frac{1}{p}\sum_{0\leq n<p}t_p(n)\varphi(f^n(x,y)).
  $$
  
  The following properties hold:
  \par
  \emph{(1)} The sequence~$(s_p)$ does not converge almost everywhere as
  $p\to +\infty$.
  \par
  \emph{(2)} If $\mathsf{P}$ is any infinite set of primes such that
  $$
  \sum_{p\in\mathsf{P}}\frac{\log p}{p}<+\infty,
  $$
  then the sequence $(s_p)_{p\in\mathsf{P}}$ converges almost everywhere
  to~$0$.
\end{proposition}

\begin{proof}
  Since $f^n(x,y)=(x+ny,y)$ for all integers $n\in\Zz$, we can compute
  $s_n$ by summing a finite geometric progression, and we obtain
  $$
  s_p(x,y)=\frac{e(x)}{p}\frac{\sin(\pi p(y-a_p/p))}{\sin(\pi
    (y-a_p/p))}
  e\Bigl(\frac{(p-1)}{2}(y-a_p/p)\Bigr).
  $$

  It follows that $s_p(x,y)\to 0$ along any infinite set $\mathsf{P}$ of
  primes such that
  $$
  \lim_{\substack{p\to +\infty\\p\in \mathsf{P}}}
  \ p\Bigl|y-\frac{a_p}{p}\Bigr|=+\infty.
  $$

  Thus~(2) follows because the assumption there implies that almost all
  $(x,y)$ satisfy
  $$
  \Bigl|y-\frac{a_p}{p}\Bigr|\geq \frac{\log p}{p}
  $$
  for all but finitely many $p\in\mathsf{P}$, by the easy
  Borel--Cantelli lemma.

  We now prove~(1). Note that if $(s_p)$ converges almost everywhere,
  the limit must be zero according to~(2), applied to a suitably sparse
  set of primes. But the formula for $s_p$ and
  the defining property of $(a_p)$ imply that for all~$x$ and almost all
  $y\in\Rr/\Zz$, we have $|s_p(x,y)|\gg 1$ for infinitely many
  primes~$p$. Thus there is almost surely a subsequence that does not
  converge to~$0$.
\end{proof}

\begin{remark}
  For more details on the type of ergodic averages considered here, we
  refer to the preprint~\cite{kowalski_ergodic}. The key fact is that
  it is proved there that the second part of the proposition
  generalizes to averages of the form
  \[
    \frac{1}{p}\sum_{0\leq n<p}t_p(n)\varphi(f^n(x))
  \]
  for suitable families $(t_p)$ of trace functions modulo~$p$ (with
  bounded conductor and without Artin--Schreier components), arbitrary
  dynamical systems $(X,\mu,f)$ with trivial Kronecker factor, and
  integrable function $\varphi\colon X\to\Cc$, \emph{under the
    proviso} that the limit as $p\to +\infty$ is done along
  sets~$\mathsf{P}$ of primes with
  $$
  \sum_{p\in\mathsf{P}}\frac{(\log p)^2}{p}<+\infty.
  $$
  
  The proposition is our current best attempt at trying to determine
  if this type of restriction on~$\mathsf{P}$ is necessary or not.
\end{remark}

\section{Preliminaries}

We summarize here some results that will be used in the proofs of the
main theorems.

First, we recall a result of Cassels, which for us essentially states
that the measure $\lambda(\mathcal{A}(\bm{a},c))$ is independent of
the value of~$c>0$.

\begin{proposition}[Cassels's Lemma]\label{CasselsLemma}
  Let~$\delta>0$ be a fixed real number and let $(I_j)_{j\geq 1}$ be a
  sequence of intervals in $\R$ such that $\lambda(I_j)\to0$ as
  $j\to\infty$. For~$j\geq 1$, let $U_j\subset I_j$ be a Lebesgue
  measurable sets such $\lambda(U_j)\geq \delta\lambda(I_j)$. We then
  have
  \[
    \lambda\Bigl(\limsup_{j\to\infty} I_j\Bigr) = \lambda\Bigl(
    \limsup_{j\to\infty} U_j\Bigr).
  \]
\end{proposition}

This is a slight modification of the result of Cassels in~\cite[Lemma
9]{C50}, and can be explicitly found in a paper of Beresnevich and
Velani~\cite[Lemma 1]{BV08}. We note that the essential content of the
proof is an application of Lebesgue's Density Theorem, and is thus
analytically non-trivial.

\begin{corollary}\label{cor-cassels}
  For any $\bm{a}\in\Omega$, the value of
  $\lambda(\mathcal{A}(\bm{a},c))$ is independent of the choice of
  $c>0$.

  In particular, for any $c > 0$, we have
  $\lambda(\mathcal{A}(\bm{a})) = \lambda(\mathcal{A}(\bm{a},c))$.
\end{corollary}

\begin{proof}
  Let~$0<c_1<c_2$ be real numbers. For any prime $p$, let
  \[
    I_p = \Bigl[\frac{a_p}{p} - \frac{c_2}{p}, \frac{a_p}{p}
    +\frac{c_2}{p}\Bigr],\quad\quad U_p = \Bigl[\frac{a_p}{p} -
    \frac{c_1}{p}, \frac{a_p}{p} +\frac{c_1}{p}\Bigr].
  \]

  We can apply Proposition~\ref{CasselsLemma} to these sets with
  $\delta = c_1/c_2$, and we obtain
  \[
    \lambda\Bigl(\limsup_{p \to \infty} I_p\Bigr) =
    \lambda\Bigl(\limsup_{p \to \infty} U_p\Bigr),
  \]
  which implies   $\lambda(\mathcal{A}(\bm{a},c_1)) =
  \lambda(\mathcal{A}(\bm{a},c_2))$
  since  
  \[
    \mathcal{A}(\bm{a},c_2) = \limsup_{p \to \infty} I_p, \quad
    \mathcal{A}(\bm{a},c_1) = \limsup_{p \to \infty} U_p,
  \]
  by definition.
  
  Since the sets $\mathcal{A}(\bm{a},c)$ are decreasing as functions
  of~$c>0$ and 
  \[
    \mathcal{A}(\bm{a})=\bigcap_{\substack{c>0\\c\in\Qq}}\mathcal{A}(\bm{a},c),
  \]
  the last assertion follows from the above.
\end{proof}

Similarly, we can use Cassels's Lemma to slightly shift the centers of
the intervals under consideration, which will be convenient in some
arguments (e.g. to remove the floor function in the definitions of the sequences in Theorem \ref{thm_badly}).

\begin{corollary}\label{cor-cassels-2}
  Let $(x_n)_{n\geq 1}$ and~$(y_n)_{n\geq 1}$ be sequences in
  $[0,1]$. Let~$(\delta_n)_{n\geq 1}$ and $(\eta_n)_{n\geq 1}$ be
  sequences of positive real numbers, with $\delta_n\to 0$ as
  $n\to+\infty$. Define
  $$
  \mathcal{X}=\limsup_n [x_n-\delta_n,x_n+\delta_n],\quad\quad
  \mathcal{Y}=\limsup_n [y_n-\eta_n,y_n+\eta_n].
  $$

  If $|x_n-y_n|\ll \delta_n$ and $\delta_n\ll \eta_n \ll \delta_n$ for
  all~$n\geq 1$, then $\lambda(\mathcal{X})=\lambda(\mathcal{Y})$.
\end{corollary}

\begin{proof}
   Since we
  have
  $$
  \lambda(\mathcal{X})=\lambda(\limsup[x_n-c\delta_n,x_n+c\delta_n])
  $$
  by Proposition~\ref{CasselsLemma} for every $c > 0$, choosing $c$ sufficiently large ensures by the assumptions made that 
  $[y_n-\eta_n,y_n+\eta_n] \subseteq [x_n-c\delta_n,x_n+c\delta_n]$. This shows
  $\lambda(\mathcal{Y}) \leq \lambda(\mathcal{X})$. The other direction now follows with exchanged roles of $\delta_n, \eta_n$ respectively $x_n,y_n$.
\end{proof}

The proof that certain sets have measure~$1$ is often obtained as a
consequence of some variant of the non-trivial direction of the
Borel--Cantelli Lemma. We will make use of the following version,
which is a very special case of a recent result of Beresnevich, Velani and the first-named author~\cite[Th.\,6]{BHV24}.

\begin{lemma}
  \label{bhv_lem}
  Let $(E_i)_{i\geq 1}$ be a sequence of Lebesgue-measurable subsets
  of $[0,1)$.  Suppose that there exist positive real numbers $C'$ and
  $\aaa$ and a sequence $(\cS_k)_{k}$ of finite subsets of~$\Zz$ such
  that the following conditions hold:
  \begin{gather}
    \label{eqn00}
    \lim_{k\to+\infty}\min\cS_k=+\infty
    \\
    \label{eqn04}
    \sum_{i\in\cS_k}\lambda(E_i) \ge \aaa,
    \\
    \label{eqn05}
    \text{ for all sufficiently large $k\geq 1$, we have}
    \\
      \sum_{\substack{s<t\\[0.5ex] s,t\in\cS_k}}
    \lambda\Bigl(E_s\cap E_t \Bigr) \leq  C'
    \Bigl(\sum_{i\in\cS_k}\lambda(E_i)\Bigr)^2\notag,
    \\
    \label{vb89}
    \text{for any $\delta>0$ and any  interval $I = [a,b] \subset
      [0,1)$, there exists $i_0$ such that }
    \\
    \lambda\Bigl(I\cap E_i\Bigr)\leq
    (1+\delta)\lambda\left(I\right)\lambda(E_i)   
    \quad
    \text{ for all $i  \geq i_0$}.\notag
  \end{gather}
  
  Then we have
  \[
    \lambda(\limsup E_i)=1.
  \]
\end{lemma}

\section{An analogy with sieve}

Let $\bm{a}\in\Omega$ and~$c>0$. For $p$ prime, let
$$
I_p=\Bigl[\frac{a_p}{p}-\frac{c}{p},\frac{a_p}{p}+\frac{c}{p}\Bigr].
$$

A real number~$x\in [0,1]$ does not
belong to $\mathcal{A}(\bm{a},c)$ if and only if, for some $Q\geq 2$,
we have
$$
x\in \bigcap_{p>Q} ([0,1]\setminus I_p).
$$

For a given~$Q$, the set
$$
\bigcap_{p>Q} ([0,1]\setminus I_p)
$$
can be viewed as a type of ``sifted set'', in analogy with the classical
sifted sets of integers, which have the form
$$
\mathcal{S}=\bigcap_{p\in\mathcal{P}}(\{1,\ldots,N\}\setminus J_p)
$$
for some set of primes~$\mathcal{P}$ (often finite) and some sets
$$
J_p=\{n\geq 1\,\mid\, n\mods{p}\in\Omega_p\}
$$
of ``excluded'' residue classes modulo~$p$. We usually wish to
show that ``most'' $x$ do belong to $\mathcal{A}(\bm{a},c)$, and this
translates to attempting to find an upper-bound for the measure of the
``sifted set'' above. Thus, we have a situation analogous to
upper-bound sieves.

Moreover, recall that a classical sieve such that $\Omega_p$ has size of
about~$\kappa$ (on average) is called a sieve of dimension~$\kappa$
(the case $\kappa=1$, the so-called linear sieve, corresponding
e.g. to primes in arithmetic progressions). In our analogy, we see
that $\kappa$ corresponds to~$2c$.

We will now give an elementary proof of the basic fact that
$\lambda(\mathcal{A}(\bm{a},c))=1$ for some~$\bm{a}$ using this analogy
and one of the simplest non-trivial sieve statements (for instance those
used by Lubotzky and Meiri in discrete groups, as
in e.g.~\cite{lubotzky-meiri}). The interesting point from the
sieve-theoretic point of view is that the usual ``sieve axioms'' are
obtained probabilistically -- on average over~$\Omega$, and not for a
fixed~$\bm{a}$.

\begin{remark}
  It would be of some interest to know if the analogy goes at all
  deeper, e.g., if there are also matching lower bounds for the measure
  of
  $$
  \bigcap_{Q<p\leq R} ([0,1]\setminus I_p)
  $$
  for suitable integers $Q$ and~$R$, for random~$\bm{a}$, or for some
  fixed sequences.
\end{remark}

\begin{proposition}\label{pr-easy}
  Let~$c>0$. There exists $\bm{a}\in\Omega$ such that
  \[
    \lambda(\mathcal{A}(\bm{a},c))=1.
  \]
\end{proposition}

The proof will not use Cassels's Lemma or any other ingredient except
for the easy part of the Borel--Cantelli lemma. Readers interested in
more precise results such as Theorem~\ref{th-random} may skip the
remainder of this section.

Fix~$c>0$.  For a sequence $\uple{a}=(a_p)$ in~$\Omega$, and for real
parameters~$X$ and~$Y$ with $1\leq X<Y$, we consider the set
$$
\mathcal{A}_{X,Y}(\uple{a},c)=\Bigl\{ x\in [0,1]\,\mid\,
\Bigl|x-\frac{a_p}{p}\Bigr|> \frac{c}{p}\text{ for all primes $p$ with $X< p\leq Y$}
\Bigr\}.
$$

The following result is a kind of ``sieve on average'' over $\Omega$. We
denote
\[
  H_{X,Y}=\sum_{X< p\leq Y}\frac{1}{p}
\]
for $1\leq X<Y$. Below, expectation and probability are taken with
respect to the natural probability measure on~$\Omega$, which was
defined above the statement of Theorem~\ref{th-random}.
 
\begin{lemma}\label{lm-sieve}
  We have
  $$
  \expect(\lambda(\mathcal{A}_{X,Y}(\bm{a},c))) \ll \frac{1}{H_{X,Y}},
  $$
  where the implied constant depends only on~$c$.
\end{lemma}

\begin{proof}
  We denote by $\varphi_p\colon [0,1]\to \{0,1\}$ the characteristic
  function of the interval $I_p(a_p)=[a_p/p-c/p,a_p/p+c/p]$, each being
  viewed as random variables on~$\Omega$ (the $\varphi_p$ are random
  functions, the $I_p$ are random intervals). Let
  $$
  N_{X,Y}=\sum_{X< p\leq Y}\varphi_p,
  $$
  again a random variable on~$\Omega$. We denote also
  $$
  \nu_{X,Y}=\int_0^1N_{X,Y}(x)dx
  $$
  and note that $\nu_{X,Y}=2cH_{X,Y}$, independently of the value
  of~$\uple{a}$.

  Noting that $\mathcal{A}_{X,Y}$ is the set of those~$x$ where
  $N_{X,Y}=0$, we deduce from Markov's inequality (on $[0,1]$ with the
  Lebesgue measure) the upper bound
  $$
  \lambda(\mathcal{A}_{X,Y})
  \leq
  \lambda\Bigl(\Bigl\{
  x\in[0,1]\,\mid\, |N_{X,Y}(x)-\nu_{X,Y}|\geq \nu_{X,Y}
  \Bigr\}\Bigr)\leq \frac{\alpha_{X,Y}}{\nu_{X,Y}^2}
  $$
  where
  $$
  \alpha_{X,Y}=\int_0^1\Bigl(N_{X,Y}(x)-\nu_{X,Y}\Bigr)^2dx
  $$
  (again a random variable on~$\Omega$).

  We have
  $$
  \alpha_{X,Y}=\int_0^1\Bigl(\sum_{X\leq
    p\leq Y}\Bigl(\varphi_p(x)-\frac{2c}{p}\Bigr)\Bigr)^2dx=
  \sum_{X< p_1,p_2\leq Y} \int_0^1
  \Bigl(\varphi_{p_1}(x)-\frac{2c}{p}\Bigr)
  \Bigl(\varphi_{p_2}(x)-\frac{2c}{p}\Bigr)
  dx.
  $$

  For $p_1=p_2$, the integral is equal to
  $$
  \int_0^1
  \Bigl(\varphi_{p_1}(x)-\frac{2c}{p_1}\Bigr)^2dx=
  \frac{2c}{p_1}\Bigl(1-\frac{2c}{p_1}\Bigr)\leq \frac{2c}{p_1}
  $$
  (variance of a Bernoulli random variable with probability of success
  $2c/p_1$), again independently of~$\uple{a}$. Thus
  $$
  \sum_{X< p_1\leq Y}\expect\Bigl(\int_0^1
  \Bigl(\varphi_{p_1}(x)-\frac{2c}{p_1}\Bigr)^2dx\Bigr)
  \leq 2cH_{X,Y}.
  $$

  We now suppose that $p_1\not=p_2$.  We then have
  \begin{equation}\label{eq-var1}
    \int_0^1 \Bigl(\varphi_{p_1}(x)-\frac{2c}{p_1}\Bigr)
    \Bigl(\varphi_{p_2}(x)-\frac{2c}{p_2}\Bigr)dx =\lambda(I_{p_1}\cap
    I_{p_2})-\frac{4c^2}{p_1p_2},
  \end{equation}
  where the first term depends on~$\uple{a}$.

  We next estimate the expectation
  $$
  \expect\Bigl( \lambda(I_{p_1}\cap I_{p_2}) \Bigr)
  $$
  over~$\uple{a}$. For this purpose, we may (and do) assume that
  $p_1<p_2$. We then have the formula
  $$
  \expect\Bigl( \lambda(I_{p_1}\cap I_{p_2}) \Bigr)= \frac{1}{p_1p_2}
  \sum_{0\leq a<p_1}\lambda\Bigl( I_{p_1}(a)\cap \bigcup_{0\leq
    b<p_2}\Bigl[ \frac{b}{p_2}-\frac{c}{p_2},
  \frac{b}{p_2}+\frac{c}{p_2} \Bigr] \Bigr).
  $$

  Drawing a picture if need be, we get
  $$
  \expect\Bigl( \lambda(I_{p_1}\cap I_{p_2}) \Bigr)=
  \frac{1}{p_1p_2}\times p_1\times \Bigl( \frac{4c^2}{p_1}+O\Bigl(
  \frac{1}{p_2} \Bigr) \Bigr)=\frac{4c^2}{p_1p_2} +O\Bigl(
  \frac{1}{p_2^2} \Bigr).
  $$

  Combined with~(\ref{eq-var1}), this leads to
  $$
  \sum_{X\leq p_1<p_2\leq Y} \expect
  \Bigl(\int_0^1 \Bigl(\varphi_{p_1}(x)-\frac{2c}{p_1}\Bigr)
  \Bigl(\varphi_{p_2}(x)-\frac{2c}{p_2}\Bigr)dx \Bigr)\ll H_{X,Y}.
  $$

  Multiplying by two to account for the case $p_1>p_2$ and adding the
  contribution where $p_1=p_2$, we conclude that
  $\expect( \alpha_{X,Y})\ll H_{X,Y}$, and hence
  $$
  \expect(\lambda(\mathcal{A}_{X,Y}))\ll H_{X,Y}^{-1},
  $$
  as claimed.
\end{proof}

We now conclude the proof of Proposition~\ref{pr-easy}.  Since
\[
  \lim_{Y\to +\infty}H_{X,Y}=\sum_{p>X} \frac{1}{p}=+\infty
\]
for any $X\geq 1$, Lemma~\ref{lm-sieve} implies that for any $X\geq 2$
and any $\eps>0$, we can find a tuple $(a_p)_{X< p\leq Y}$, with
$0\leq a_p<p$,  such that
\[
  \lambda\Bigl(\Bigl\{ x\in [0,1]\,\mid\,
  \Bigl|x-\frac{a_p}{p}\Bigr|\leq \frac{c}{p}\text{ for some prime $p$
    with $X< p\leq Y$}
  \Bigr\}\Bigr)\geq 1-\eps.
\]

Let $X_1=1$. Apply the previous remark first with (say) $X=1$ and
$\eps=1/2$, and denote $X_2$ a suitable value of~$Y$. Then apply the
assumption with $X=X_2$ and $\eps=1/4$, calling~$X_3$ the value
of~$Y$; repeating, we obtain a strictly increasing sequence
$(X_n)_{n\geq 1}$ of integers and a sequence $(a_p)\in\Omega$ such
that the set
\[
  B_n=\Bigl\{ x\in [0,1]\,\mid\,
  \Bigl|x-\frac{a_p}{p}\Bigr|\leq \frac{c}{p}\text{ for some prime $p$
    with $X_n<p\leq X_{n+1}$} \Bigr\}
\]
satisfies $\lambda(B_n)\geq 1-2^{-n}$ for any~$n\geq 1$.

If $x\in [0,1]$ belongs to infinitely sets $B_n$, then
$x\in \mathcal{A}(\uple{a},c)$. On the other hand, since
$$
\sum_{n\geq 1}\lambda([0,1]\setminus B_n)<+\infty,
$$
the easy Borel--Cantelli Lemma shows that almost every $x\in [0,1]$
belongs at most to finitely many sets~$[0,1]\setminus B_n$. Thus we
have proved indeed that $\lambda(\mathcal{A}(\bm{a}),c)=1$ for
some~$\bm{a}\in\Omega$.

\section{Probabilistic results}

This section is devoted to the proof of
Theorem~\ref{th-random}.
We will in fact provide a slightly more precise error term than stated
in~(\ref{eq-count}).  The key tool is the following result, whose basic
idea goes back to Weyl; we refer to the account in Harman's
book~\cite[Lemma\,1.5]{Harman_1998} for more general versions and
further references.

\begin{lemma}\label{harman}
  Let $(X,\mu)$ be a probability space.  Let $(f_k)_{k\geq 1}$ be a
  sequence of non-negative $\mu$-integrable functions. Let $\nu_k=\int_{X}
  f_kd\mu$ for $k\geq 1$ and
  define
  \[
    \Psi(N) = \sum_{k = 1}^{N} \nu_k
  \]
  for $N\geq 1$.  Assume that $\Psi(N)\to+\infty$ as $N\to +\infty$, and
  that there exists a constant $c\geq 0$ such that, for integers~$M$
  and~$N$ with $1\leq M<N$, we have
  \begin{equation}\label{local_variance}
    \int_{X} \Bigl(\sum_{M \leq k < N} (f_k(x) - \nu_k)\Bigr)^2
    d\mu(x) \leq c\sum_{M \leq k < N}\nu_k.
  \end{equation}

  Let $\delta > 0$. For $\mu$-almost all $x \in X$, we have then
  \[
    \sum_{k = 1}^N f_k(x) = \sum_{k=1}^N\nu_k + O\Bigl( \Psi(N)^{1/2}
    \Bigl(\log(\Psi(N) +2)\Bigr)^{3/2 + \delta} + \max_{1 \leq k \leq
      N}\nu_k\Bigr)
  \]
  for~$N\geq 1$.
\end{lemma}

We also need the following lemma that is considered standard in the context of metric Diophantine approximation (see e.g. \cite[Chapter 4]{Harman_1998}). Because of the short proof, we provide it for the sake of completeness.

\begin{lemma}\label{coprime_overlap}
  Let $p<q$ be prime numbers and let $0 < c <\frac{1}{2}$. Then we have
  \begin{equation}\label{prime_independence}
    \int_{0}^1 \mathds{1}_{[\lVert p\alpha \rVert < c]}
    \mathds{1}_{[\lVert q\alpha \rVert < c]}
    d\alpha \leq 4c^2 + O\Bigl(\frac{c}{q}\Bigr).
  \end{equation}
\end{lemma}

\begin{proof}
  Denoting $ \delta = \frac{c}{q}$ and $\Delta = \frac{c}{p}$, we
  define a piecewise linear function~$w$ (see \cite{ABH23}) by
  \begin{equation}\label{weight_fct}
    w(y) = \begin{cases}
      2\delta &\text{ if }
      0 \leq \lvert y \rvert \leq \Delta - \delta,\\
      \Delta + \delta - y &\text{ if }
      \Delta - \delta < \lvert y \rvert \leq \Delta + \delta,\\
      0 &\text{ otherwise.}
    \end{cases}
  \end{equation}

  The function $w$ is monotonically increasing on $]-\infty,0]$ and
  monotonically decreasing on $[0,+\infty[$. Using the Chinese Remainder Theorem, we obtain
  \[
    \int_{0}^1 \mathds{1}_{[\lVert p\alpha \rVert <
      c]}\mathds{1}_{[\lVert q\alpha \rVert < c]}d\alpha \leq \sum_{x
      \in \mathbf{Z}} w\left(\frac{x}{pq}\right) \leq w(0) + \int_{\R}
    w(x)dx.
  \]
  
  Observing finally that
  \[
    \int_{\Rr} w(x)dx = 4c^2, \quad\quad w(0) = \frac{2c}{q},
  \]
 property \eqref{prime_independence} follows.
\end{proof}
  
\begin{proof}[Proof of Theorem~\ref{th-random}]
  According to Corollary~\ref{cor-cassels}, we may assume that~$c>0$
  is fixed.
  
  Let~$\Omega'$ be the probability space ${[0,1] \times \Omega}$ with the
  product measure~$\mu$ of the probability measure on~$\Omega$ and the
  Lebesgue measure.
  We define random variables $\varphi_p$ (for $p$ prime) and $N_X$
  (for $X\geq 2$) on~$\Omega'$ by
  \begin{gather*}
    \varphi_p(\alpha,\bm{a})
    = \mathds{1}_{\left[\left\lvert \alpha - \tfrac{a_p}{p} \right\rvert
        \leq \frac{c}{p}\right]},
    \\
    N_X(\alpha,\bm{a}) = \sum_{p \leq X}\varphi_p(\alpha,\bm{a}) =
    \Bigl|\Bigl\{ p \leq X\,\mid\, \left\lvert \alpha - \frac{a}{p}
    \right\rvert \leq \frac{c}{p}\Bigr\}\Bigr|.
  \end{gather*}

  We will apply Lemma~\ref{harman} to the probability space
  $(\Omega',\mu)$ and to the functions $(\varphi_p)$, which have integral
  $2c/p$.  In particular, with the notation of the lemma, we have
  $\Psi(N)\to+\infty$ as $N\to +\infty$.
  
  We begin to check \eqref{local_variance}. Let $X \leq Y$ be
  fixed integers.  By Fubini's Theorem, we can write
  \begin{align*}
    \expect\Bigl(\Bigl(\sum_{X<p\leq Y}\varphi_p\Bigr)^2\Bigr)
    =\sum_{X< p,q\leq Y}\expect(\varphi_p\varphi_q)
     &=\sum_{X< p,q\leq Y} \frac{1}{pq}\sum_{{\substack{0\leq a\leq p-1 \\0\leq
     b    \leq q-1}}} \int_0^1 \mathds{1}_{\left[\left\lvert \alpha -
     \frac{a}{p} \right\rvert \leq \frac{c}{p}\right]}
     \mathds{1}_{\left[\left\lvert \alpha - \frac{b}{q} \right\rvert \leq
     \frac{c}{q}\right]}d\alpha
     \\    &=\sum_{X< p,q\leq Y} {\frac{1}{pq}}\int_0^1
         \mathds{1}_{[\lVert p\alpha \rVert < c]}\mathds{1}_{[\lVert q\alpha
       \rVert < c]}d\alpha.
  \end{align*}

  By Lemma \ref{coprime_overlap}
  we have for $p < q$ that
  \[\int_0^1
         \mathds{1}_{[\lVert p\alpha \rVert < c]}\mathds{1}_{[\lVert q\alpha
       \rVert < c]}d\alpha \leq 4c^2 + O(c/q)\]
  
  Thus we obtain
  \begin{align*}
    \expect\Bigl(\Bigl(\sum_{X<p\leq Y}\varphi_p\Bigr)^2\Bigr)
    & = 2\sum_{X< p < q \leq Y}\frac{4c^2 + O(c/q)}{pq} + \sum_{X < p \leq Y}\expect(\varphi_p^2) 
    \\&\leq \Bigl(\sum_{X< p \leq Y}\frac{2c}{p}\Bigr)^2 + 
    O\Bigl(\Bigl(\sum_{X< q \leq Y}\frac{1}{q^2}\Bigr)\cdot \Bigl(\sum_{X < p \leq Y}\frac{1}{p}\Bigr)\Bigr) +
    2c \sum_{X< p
      \leq Y}\frac{1}{p}
       \\
    &= \Bigl(\sum_{X< p \leq Y}\expect(\varphi_p)\Bigr)^2 + O\Bigl(\sum_{X<
      p \leq Y}\expect(\varphi_p)\Bigr).
  \end{align*}

  Thus Lemma~\ref{harman} implies that there exists a set
  $A \subset \Omega'$ with $\mu(A) = 1$ such that
  \begin{equation}\label{asmyp_psi}
    N_X(\alpha,\bm{a}) = 
    \Psi(X) +
    O\Bigl(\sqrt{\Psi(X)}\Bigl(\log (\Psi(X))^{3/2 + \delta}\Bigr)\Bigr)
  \end{equation}
  for all $(\alpha,\bm{a}) \in A$, where
  \[
    \Psi(X) = 2c \sum_{p \leq X}\frac{1}{p}.
  \]
  
  We finally deduce Theorem~\ref{th-random} by a standard argument: by
  Fubini's theorem, the function~$F$ on~$\Omega$ defined by
  \[
    F(\bm{a}) = \int_{0}^1 \mathds{1}_A(\alpha,\bm{a})d\alpha
  \]
  for $\bm{a}\in\Omega$ is defined almost everywhere and satisfies
  $0 \leq F \leq 1$.  Since $\expect(F)=1$ (where expectation is
  computed on $\Omega$), this implies $F=1$ almost everywhere. This
  means that for almost every $\bm{a}\in\Omega$, there exists a set
  $S(\bm{a})\subset [0,1]$ of Lebesgue measure~$1$ such
  that~\eqref{asmyp_psi} holds for any $\alpha \in S$. By the Mertens
  formula, we have
  \[
    \Psi(X) = 2c \log \log X + O(1),
  \]
  and the theorem follows.
\end{proof}

\begin{proof}[Proof of Corollary \ref{uniform_c}]
  We use the notation from the previous proof, except that we write
  \[
    N_X^{(c)}(\alpha,\bm{a})
    = \Bigl|\Bigl\{ p \leq X\,\mid\, \left\lvert \alpha - \frac{a_p}{p}
    \right\rvert \leq \frac{c}{p}\Bigr\}\Bigr|
  \]
  since the dependency on~$c$ will play a role.
  
  Since $\mathbf{Q}$ is countable, an application of
  Theorem~\ref{th-random} shows that for almost every sequence
  $\bm{a}\in\Omega$, there exists a set $S(\bm{a}) \subseteq [0,1]$ of
  Lebesgue measure~$1$ such that for any $\alpha \in S(\bm{a})$, any
  $\delta > 0$ and any $c \in [0,1/2)\cap \mathbf{Q}$, we have
  \[
    \Bigl|\Bigl\{ p \leq X\,\mid\, \left\lvert \alpha - \frac{a_p}{p}
    \right\rvert \leq \frac{c}{p}\Bigr\}\Bigr| = 2c\log \log X +
    O\Bigl(\sqrt{\log \log X}\Bigl(\log\log \log(X)\Bigr)^{3/2 +
      \delta}\Bigr)
  \]
  as $X \to \infty$.

  For irrational $c$ and any $\eps>0$, pick $c_1 < c < c_2$ such
  that $c_1,c_2 \in \mathbf{Q}$ and
  $\frac{c_2}{c_1} \leq 1 + \varepsilon$. Since
  $N_X^{(c)}(\alpha,\bm{a})$ is monotonically increasing in $c$, we have
  \[
    \frac{N_X^{(c_1)}(\alpha,\bm{a})}{2c\log \log X} \leq
    \frac{N_X^{(c)}(\alpha,\bm{a})}{2c\log \log X} \leq
    \frac{N_X^{(c_2)}(\alpha,\bm{a})}{2c\log \log X},
  \]
  which implies that
  \[1-\varepsilon \leq \liminf_{X \to \infty}
    \frac{N_X^{(c)}(\alpha,\bm{a})}{2c\log \log X} \leq \limsup_{X \to
      \infty} \frac{N_X^{(c)}(\alpha,\bm{a})}{2c\log \log X} \leq 1 +
    \varepsilon.
  \]
  
  The corollary follows by letting $\eps \to 0$.
\end{proof}

\section{The greedy sequence}\label{sec-greedy}

We will explicitly construct a sequence $\bm{g}=(a_p)_{p}$ such that
$\mathcal{A}(\bm{g},2) = [0,1[$. By Corollary~\ref{cor-cassels}, this
implies that $\lambda\left(\mathcal{A}(\bm{g},c)\right) = 1$ for any
$c > 0$, which then proves the first part of
Theorem~\ref{thm_greedy}. The final statement of Theorem~\ref{thm_greedy},
concerning the Hausdorff dimension of the exceptional set, is left for
Section~\ref{sec-HD}.

The idea is to use a ``greedy'' construction, iterating over positive
integers~$j$ so that each iteration step constructs $a_p$ for a
range~$\mathcal{I}_j$ of primes~$p$ in such a way that
\[
  \bigcup_{p\in\mathcal{I}_j}\Bigl]
  {\frac{a_p}{p}-\frac{2}{p},\frac{a_p}{p}+\frac{2}{p}}
  \Bigr[ \supseteq [0,1[,
\]
and the ranges are disjoint at each step (if~$j_1<j_2$, then we have
$p<q$ for any primes $p\in\mathcal{I}_{j_1}$
and~$q\in\mathcal{I}_{j_2}$).  The iteration is also performed so that
we cover the interval from left to right.

We formally describe an algorithm that performs this construction
starting from a given prime~$p_0$.
For any real number $x$, we will denote by $\nextp(x)$ the smallest
prime number strictly larger than $x$.

\begin{algorithm}
  \caption{Greedy algorithm: cover the interval once starting
    from~$p_0$}
  \begin{algorithmic}
    \Procedure{CoverOnceFrom}{$p_0$}
    \State $q \gets p_0-1$
    \State $x \gets 0$
    \While{$x < 1$}
    \State $q \gets \nextp(q)$
    \State $a_q := \lfloor qx \rfloor$
    \State $x \gets \frac{a_q}{q}+\frac{2}{q}$
    \EndWhile
    \State \textbf{return} $(a_{p_0},\ldots,a_q)$
    \EndProcedure
  \end{algorithmic}
\end{algorithm}

The following lemma gives the required analysis of this algorithm and
concludes the proof that~$\mathcal{A}(\bm{g},2)=[0,1[$.

\begin{lemma}\label{lem_greedy}
  For any prime~$p_0$, the above algorithm terminates. Furthermore, if
  $(a_{p_0},\ldots,a_q)$ are the integers returned by the procedure,
  then we have
  \[
    \bigcup_{p_0\leq p\leq q}\Bigl] {\frac{a_p}{p}-\frac{2}{p},\frac{a_p}{p}+\frac{2}{p}
    \Bigr[=[0,1].}
  \]
\end{lemma}

\begin{proof}
  For any prime~$p$ with~$p_0\leq p$ for which
  $a_p$ has been defined by running the algorithm, denote
  \[
    A(p)=\Bigl] {\frac{a_p}{p}-\frac{2}{p},\frac{a_p}{p}+\frac{2}{p} } \Bigr[.
  \]

  We note that~$0\in A(p_0)$ by construction.

  Let~$q$ and~$x$ be the values of the variables after the first step
  of some iteration of the \textbf{while} loop (so that~$q$ is prime)
  and~$q'$, $x'$ the values at the same point of the next iteration,
  if there is one.
  We claim that
  \begin{enumerate}
  \item[(i)] $x' \geq x + \frac{1}{q'}$,
  \item[(ii)] $A(q') \cap A(q) \neq \emptyset$.
  \end{enumerate}

  To prove (i), we observe simply that 
  \[
    x' = \frac{\lfloor q'x\rfloor + 2}{q'} \geq
    \frac{q'x +1}{q'} = x+ \frac{1}{q'},
  \]
  and (ii) follows from
  \[
    x'=
    {\frac{\lfloor q'x\rfloor +2}{q'}\leq \frac{q'x+2}{q'} \leq x+\frac{2}{q}.}
  \]

  Since we have
  \[
    \sum_{p\geq p_0}\frac{1}{p}=+\infty,
  \]
  the first claim implies that the \textbf{while} loop indeed
  terminates after finitely many steps.  The second implies that the union
  \[
    \bigcup_{p_0\leq p\leq q}A(p)
  \]
  is connected, hence is an interval. Since it contains~$0\in A(p_0)$
  and~$1\in A(q)$, it is equal to~$[0,1]$, as desired.
\end{proof}

\begin{remark}
  Here are some data on the results of running this algorithm.  The
  five first coverings of the interval are obtained for the following
  ranges of primes:
  \begin{gather*}
    \{2\},\quad\quad \{3,5\},\quad\quad  \{7,\ldots,47\},\\
    \{53,\ldots, 2593\},\quad\quad
    \{2609,\ldots, 4835851\}.
  \end{gather*}

  The pairs~$(p,a_p)$ for primes $p\leq 47$ are as follows:
  \[
    (7, 0 ),(11, 3 ),(13, 5 ),(17, 9 ),(19, 12 ),(23, 16),(29 22),
    (31, 25),(37, 32),(41, 37),(43, 40),(47, 45).
  \]
\end{remark}

\section{Deterministic sequences}
\label{sec-deterministic}

In this Section, we prove Theorem \ref{th-deterministic}.  We will
first observe that the case of~$\bm{a}(\beta)$ can be handled using
results of twisted diophantine approximation. We then consider the
case of $\bm{b}(\beta)$, where no corresponding result appears to be known.

\subsection{The rotation case}

We first consider the sequence $\bm{a}(\beta)=(a_p)$, where we recall
that~$a_p=\lfloor \beta pn_p\rfloor\pmod{p}$.  It turns out that there
is a connection between the problem in the case of $\bm{a}(\beta)$ and
the area of \emph{twisted diophantine approximation}, from which we can
quote to obtain the result in that case (one could also treat this case
in the same manner that we will handle the sequence $(b_p)$ {in} the next
section).

In twisted diophantine approximation, one considers a decreasing
function~$\psi$ defined on positive integers and a real
number~$\beta$, and the question is to understand the properties of
the set
\[
  A_{\beta}(\psi) = \left\{x \in [0,1): \lVert n\beta - x \rVert \leq
    \psi(n) \text{ for infinitely many integers } n\geq 1\right\}.
\]

{ 
Using the Prime
Number Theorem (in the very weak form of $c_1 n \log n \leq p_n \leq c_2 n \log n$ for some constants $c_1,c_2 > 0$) and setting $\psi(1)=1$ and $\psi(n)=1/(n\log n)$
for~$n\geq 2$, we observe that
}
{
\[\frac{a_p}{p} = \{n_p\beta\} + O(1/p) =  \{n_p\beta\} + O(\psi(n_p)).\]
}
{Thus by Cassels's Lemma in the form of Corollary~\ref{cor-cassels-2}, we have
\begin{equation}\label{pnt_cass}
  \lambda(A_{\beta}(\psi)) =
  \lambda(\mathcal{A}(\bm{a}(\beta))).
\end{equation}
Since $\psi$ is monotonically decreasing and $\sum_{n \in \mathbb{N}}\psi(n) = \infty$, the result of Kurzweil~\cite[Th.\,1]{K55} implies that $ \lambda(A_{\beta}(\psi)) = 1$ (and therefore by \eqref{pnt_cass}
$\lambda(\mathcal{A}(\bm{a}(\beta))) = 1$) if $\beta$ is badly
approximable.}

{
But more generally, a theorem of Fuchs and Kim (see~\cite[Theorem 2]{FK16}) states that $\lambda(A_\beta(\psi)) = 1$ (and in view of \eqref{pnt_cass} equivalently $\lambda(\mathcal{A}(\bm{a}(\beta))) = 1$)}
{if and only if}~$\beta$ satisfies the condition
\[
  \sum_{k\geq 0}\sum_{q_k\leq
    n<q_{k+1}}\min(\psi(n),\|q_k\beta\|)=+\infty.
\]

The above property is indeed satisfied for $\psi(n) = \frac{1}{n \log n}$ and $\beta$ i.b.a.:

{Since $\lVert q_k\beta\rVert \geq \frac{1}{a_{k+1} + 2}\frac{1}{q_k}$, it follows that 
for $a_{k+1} \leq B$, we have for sufficiently large $k$, $\lVert q_k\beta\rVert \geq \frac{1}{n \log n}$ for all $q_k \leq n \leq q_{k+1} \leq (B+1)q_k$. 
Let $(j_k)_{k \geq 1},(T_k)_{k \geq 1}$ being the sequences along which 
\eqref{bad_range} and \eqref{range_length} hold. We may assume without loss of generality that for all $k \geq 1$, $j_k + T_k < j_{k+1}$ since otherwise we simply pass to a sufficiently sparse subsequence. Further, using the fact that $(q_j)_{j \in \mathbb{N}}$ grows at least exponentially, we have
\[q_{j_k+T_k} \geq q_j^{1+\delta'}, \quad k \geq 1\]
for some $\delta' = \delta'(\delta) > 0$.
Therefore, we obtain
}

\[\begin{split}  \sum_{k\geq 0}\sum_{q_k\leq
    n<q_{k+1}}\min\left(\frac{1}{n \log n},\|q_k\beta\|\right)
    &\geq 
    \sum_{k\geq 0}\sum_{i = 0}^{T_k}\sum_{q_{k+i}\leq
    n<q_{k+i+1}} \frac{1}{n \log n}
    \\&\geq  \sum_{k\geq 0} \sum_{q_{j_k} \leq n \leq q_{j_k}^{1 + \delta'}}\frac{1}{n \log n} \gg_{\delta'} \sum_{k\geq 0} 1 = \infty.
    \end{split}\]

This concludes the proof for~$\bm{a}(\beta)$.

\subsection{The prime rotation case}

We now consider the case of~$b_p=\lfloor \beta
p^2\rfloor\pmod{p}$. One can also connect this to twisted diophantine
approximation, as in the previous section, but there is no established
result to quote here. {Thus we will establish Theorem~\ref{prime_twisted} ourselves, and will then use the link as before to prove Theorem~\ref{th-deterministic}. This will be done in Subsection~\ref{ssec-prime-twisted}.} {We now
begin the argument towards the proof of Theorem~\ref{prime_twisted}, with} the properties of the primes entering through two fundamental
properties, which we now state.

\begin{proposition}\label{pr-primes}
  \begin{enumth}
  \item Let~$\beta$ be an irrational number. The fractional parts
    $\{p\beta\}$ for $x<p\leq 2x$ are equidistributed in $[0,1]$ for
    the Lebesgue measure as $x\to+\infty$.
  \item For any integer~$h\geq 1$ and real number $x\geq 2$ with
    $h\leq x^2$, we have
    \[
      |\{p\leq x,\mid\, p+h\text{ is prime }\} | \ll
      \frac{h}{\varphi(h)}\frac{x}{(\log x)^2}
    \]
    where the implied constant is absolute.
  \end{enumth}
\end{proposition}

The first of these results is due to Vinogradov (see for
instance~\cite[Th.\,21.3]{ant} for a proof); the second is a standard
application of basic (upper-bound) sieve theory (see for
instance~\cite[Th.\,6.7]{ant}).

\begin{definition}[Badly approximable range]
  Let $\beta = [0;a_1,a_2,\ldots]$ be an irrational number and
  $B \geq 1$ an integer. For integers $U<V$, we say that the interval
  $[U,V]$ is a \textit{$B$-badly approximable range} if we have for
  $a_{t} < 2^U < a_{t+1} < \cdots < a_{t+r-1} < 2^V < a_{t+r}$ that
  $\max_{0 \leq j \leq r}a_{t + j} \leq B$.
\end{definition}

Let~$\beta$ be an i.b.a. real number, and
let~$B \geq 2$ and~$\delta \in ]0,1[$ be the parameters appearing in the
definition. The estimates below will depend on these values, but
since~$\beta$ is fixed, we will treat them as constants.

It follows straightforwardly from the definition of $\beta$ being i.b.a. that we can find a
sequence of intervals $([U_k,V_k])_{k\geq 1}$ such that
\begin{align*}
  &\text{The sequences $(U_k)_{k \in \N}$ and $(V_k)_{k \in \N}$ are
    monotonically increasing,}
  \\
  &[U_k,V_k] \cap [U_j,V_j] = \emptyset \text{ for
    all } k \neq j,
  \\
  &[U_k-2,V_k+2] \text{ is a $B$-badly approximable range for all~$k$,}
  \\
  &(2^{U_k})^{1+\delta/2}< 2^{V_k}< (2^{U_k})^2
    \text{ for all
    $k$.}
\end{align*}

We note that the last condition implies in particular
\[
  \frac{V_k}{U_k}\geq 1+\frac{\delta}{2}.
\]

   {For $c > 0$, let 
\[
	\mathcal{A}(\beta,c) := \{x \in [0,1)\,\mid\, \liminf_{p \to \infty} p\lVert
	p\beta - x \rVert < c\bigr\}.
	\]
    }

As before, in view of Corollary \ref{cor-cassels}, it is enough to prove that $\lambda(\mathcal{A}(\beta,c))=1$
for some~$c>0$. We pick now an arbitrary value $c$ such that
\[
0<c<\frac{1}{10B}.
\]

We denote by~$\mathcal{I}$ the set of integers defined by
\[
  \mathcal{I} = \bigcup_{k \geq 1}[U_k,V_k].
\]

For~$p$ prime and $i \in \mathcal{I}$, we define
\[
  A_{p,i} = \Bigl[p\beta - \frac{c}{2^{i+1}}, p\beta +
    \frac{c}{2^{i+1}}\Bigr],
\]
and we put
\[
  E_i= \bigcup_{2^i \leq p < 2^{i+1}}A_{p,i}.
\]

For completeness, we define $E_i = \emptyset$ when
$i \notin \mathcal{I}$.

It is immediate that $\mathcal{A}(\beta,c)$ contains
$\limsup E_i$, so it suffices to prove that
\[
  \lambda(\limsup_{i}E_i)=1.
\]

We prove this using Lemma~\ref{bhv_lem}, which we apply to the
sequence of subsets $(E_i)$ and the sets of integers
$\cS_k=[U_k,V_k]$.

The condition~(\ref{eqn00}) is immediate. We next check the conditions
\eqref{eqn04} and \eqref{vb89}.

\begin{lemma}\label{asymp_ind}
  \begin{enumth}
  \item For any $i \in \mathcal{I}$, we have
    \begin{equation}\label{disj_union}
     { \lambda(E_i) = \sum_{2^i \leq p < 2^{i+1}}\lambda(A_{p,i}).}
    \end{equation}
  \item For $k$ sufficiently large, we have
    \begin{equation}
      \label{accum_meas}
      \sum_{i \in [U_k,V_k]}\lambda(E_i) \gg 1,
    \end{equation}
    where the implied constant depends on~$\delta$.
  \item For any interval~$I$, any $\eta>0$, we have
  \begin{equation}\label{equidistr}
    \lambda\left(I\cap E_i\right)
    \leq (1+\eta)\lambda\left(I\right)\lambda(E_i)
  \end{equation}
  for $i$ large enough.
\end{enumth}
\end{lemma}

\begin{proof}
  To prove~(\ref{disj_union}), we show that the sets $A_{p,i}$ whose
  union gives~$E_i$ are disjoint.
  Indeed, assume that ${2^i \leq p< p'< 2^{i+1}}$ and that $A_{p,i}\cap A_{p',i}$
  is not empty.
  Ths implies by definition that
  \[
    \lVert (p'-p)\beta\rVert \leq {\frac{c}{2^i}<\frac{1}{5 B2^{i+1}}}.
  \]

  Let~$l$ be the integer such that ${q_{l}\leq 2^{i+1} <q_{l+1}}$, where as
  before we denote by $(p_l/q_l)_l$ the sequence of convergents of the
  continued fraction expansion of~$\beta$.  Since {$0<p'-p<2^{i+1}$}, the best
  approximation property of convergents implies
  \[
    \lVert (p'-p)\beta\rVert \geq \lVert q_l\beta\rVert.
  \]

  Since~$i$ belongs to some $B$-badly approximable range, on the other
  hand, we have by the well-known inequality $q_l\lVert q_l\beta\rVert \leq \frac{1}{a_{l+1}+2}$ that
  \[
    \lVert q_{l}\beta\rVert \geq \frac{1}{2Bq_l} \geq \frac{1}{{2B2^{i+1}}},
  \]
  and combining these inequalities gives a contradiction. 

  The bound~(\ref{accum_meas}) follows from~(\ref{disj_union}) and the
  Chebychev estimate:
  \[
    {\sum_{U_k\leq i\leq V_k} \lambda(E_i) \gg \sum_{U_k\leq i\leq
      V_k}\frac{1}{2^i}\frac{2^i}{i} \gg
    \log\Bigl(\frac{V_k}{U_k}\Bigr)\gg \delta\gg 1.}
  \]

  We now prove~\eqref{equidistr}. Fix an interval $I = [a,b]$
  and~$\eta>0$. Note that~(\ref{equidistr}) is trivially true for
  any~$i\notin\mathcal{I}$, so we assume that~$i\in\mathcal{I}$.
  For ~$i\in\mathcal{I}$, we have
  \[
    \lambda([a,b] \cap E_i) \leq \frac{1}{c{2^{i}}}\Bigl| \Bigl\{ p \,\mid\,
    {2^i \leq p<2^{i+1}\text{ and } \{p \beta\} \in [
    a-\tfrac{c}{2^{i+1}},b+\tfrac{c}{2^{i+1}}]} \Bigr\}\Bigr|.
  \]

  Since $\beta$ is irrational, the fractional parts of $p\beta$ become
  equidistributed in $[0,1]$ by Vinogradov's Theorem
  (Proposition~\ref{pr-primes}, (1)). Thus we have
  \[
    \Bigl| \Bigl\{ p \,\mid\, {2^i\leq p < 2^{i+1}\text{ and } \{p \beta\} \in
    [ a-\tfrac{c}{2^{i+1}},b+\tfrac{c}{2^{i+1}}] \Bigr\}\Bigr| \sim
    (b-a)
    \sum_{2^i\leq p<2^{i+1}}1}
  \]
  as $i\to+\infty$, and this implies the claim since by~(\ref{disj_union})
  \[
    \lambda(E_i)={\frac{1}{c2^{i}} \sum_{2^i\leq p \leq 2^{i+1}}1}.
  \]
\end{proof}

There only remains to check the condition \eqref{eqn05}.  This
requires some preliminaries.

For integers~$i$ and~$j$, we denote
\begin{equation}\label{bohr_def}
  \mathcal{N}(i,j)
  = \Bigl\{N\geq 1\,\mid\, N \leq 2^{j}\text{ and } \lVert N \beta
  \rVert \leq \frac{1}{2^i}\Bigr\},
\end{equation}
which are examples of \emph{Bohr sets}.  We will need to know that the
function $n/\varphi(n)$ is essentially bounded on average over such
sets. This follows intuitively from the well-known fact that Bohr sets
have additive structure, whereas the integers where~$\varphi(n)$ is
small satisfy strong multiplicative constraints. To prove the result,
we will use the relation between Bohr sets and generalized arithmetic
progressions. 

Thus, for integers $x$, $y$, $z$, we denote
\[
  P(x,y,z) = \{ax + by\,\mid\, a,b \in \Zz\text{ and } \lvert
  a\rvert,\lvert b\rvert \leq z\},
\]
which is a (rank~$2$) finite generalized arithmetic progression. The
average of $n/\varphi(n)$ over such a progression can be estimated
elementarily (this is also found in work of Chow~\cite{chow}, in a
related setup).
	
\begin{lemma}[Average over $\varphi$]\label{av_phi}
  Let $P(x,y,z)$ be a generalized arithmetic progression with $x$, $y$
  coprime and $x<y$. We have
  \[
    \sum_{n \in P(x,y,z)}\frac{n}{\varphi(n)} \ll z^2 + z\log (zy).
  \]
\end{lemma}

\begin{proof}
  Using the formula
  \[
    \frac{n}{\varphi(n)}= \sum_{d \mid n}\frac{\mu^2(d)}{\varphi(d)},
  \]
  we have
  \[
    \sum_{n \in P(x,y,z)}\frac{n}{\varphi(n)} = \sum_{d \leq zy}
    \frac{\mu^2(d)}{d}\sum_{\substack{n \in P(x,y,z)\\d \mid n}}1.
  \]
  
  For any squarefree integer~$d$, we write
  \[
    d_1 = (x,d),\quad\quad d_2 = (y,d),\quad\quad d_3 = d/(d_1d_2).
  \]

  Since~$x$ and~$y$ are coprime and~$d$ is squarefree, we have
  $d=d_1d_2d_3$, and $(d_1,d_2,d_3)$ are pairwise coprime.
  Thus the number of $n\in P(x,y,z)$ divisible by~$d$ is the number of
  integer solutions $(a,b)$, with $|a|\leq z$, $|b|\leq z$, of the
  system
  \begin{align}
    \label{d1} ax + by &\equiv 0 \pmod{d_1},\\
    \label{d2} ax + by &\equiv 0 \pmod{d_2},\\
    \label{d3} ax + by &\equiv 0 \pmod{d_3}.
  \end{align}

  Since $d_1\mid x$ and $(x,y) = 1$, \eqref{d1} is equivalent to
  $b \equiv 0 \pmod{d_1}$.  Similarly, \eqref{d2} is equivalent to
  $a \equiv 0 \pmod{d_2}$. Since $(x,d_3) = 1$, the number of
  solutions is equal to
  \[
    \sum_{\substack{\lvert b \rvert \leq z\\ b \equiv 0 \pmod{d_1}}}
    \sum_{\substack{\lvert a \rvert \leq z\\ a \equiv 0 \pmod{d_2}
        \\
        a \equiv - byx^{-1}\pmod{d_3} }} 1.
  \]

  Regardless of the value of $- byx^{-1}$, the coprimality of $d_2$
  and~$d_3$ implies that
  \[
    \sum_{\substack{\lvert a \rvert \leq z\\ a \equiv 0 \pmod{d_2}\\
        a  \equiv - byx^{-1}\pmod{d_3}
      }} 1
    \ll \frac{z}{d_2d_3} + 1,
  \]
  and arguing in the same way for the outer sum, we obtain
  \[
    \sum_{\substack{n \in P(x,y,z)\\d \mid n}}1 \ll
    \frac{z^2}{d_1d_2d_3} + \frac{z}{d_1} + \frac{z}{d_2d_3} +
      1 {\ll} \frac{z^2}{d} + z.
  \]
  
  We then obtain
  \[
    \sum_{n \in P(x,y,z)}\frac{n}{\varphi(n)} \ll z^2\sum_{d \leq yz}
    \frac{1}{d^2} + z \sum_{d \leq zy}
    \frac{1}{d},
  \]
  as claimed.
\end{proof}

The relation between Bohr sets and generalized arithmetic progressions
is given by the following lemma.

 \begin{lemma}[Structure of the Bohr sets]\label{bohr_AP}
  Assume $i\leq j$ are integers such that $[i-2,j+2]$ is a $B$-badly
  approxmable range.  There exist integers $(x,y,z)$, with $x$ and~$y$
  coprime and
  \[
    x < y \ll 2^{(j+i)/2},\quad\quad z \ll 2^{(j-i)/2},
  \]
  where the implied constants depend only on $B$, such that
  \begin{equation}
    \label{bohr_inclusion}
    \mathcal{N}(i,j) \subseteq P(x,y,z).
  \end{equation}
 
  In particular, it follows that
  \begin{equation}\label{cardinality_bohr}
    |\mathcal{N}(i,j)|\ll 2^{j-i}
  \end{equation}
  and
  \begin{equation}\label{average_bohr}
    \sum_{N\in\mathcal{N}(i,j)}\frac{N}{\varphi(N)}
    \ll     2^{j-i}.
  \end{equation}
\end{lemma}

\begin{proof}
We remark that
  \[
    \mathcal{N}(i,j) \subseteq P(q_{r},q_r + q_{r-1},z)
  \]
  for any $r\geq 1$, where
  \[
    z= \max\Bigl\{ 2^j \lVert q_r \beta \rVert + \frac{q_r}{2^i}, 2^j
    \lVert (q_r + q_{r-1})\beta \rVert + \frac{q_r + q_{r-1}}{2^i}
    \Bigr\}.
  \]
  While this fact is relatively standard, the only reference we are aware of
  is in a blog post of T. Tao ~\cite[discussion
  around~(4)]{tao_blog}.

  Choosing $r$ now such that $q_{r-1} < 2^{(j+i)/2} \leq q_r$ and
  using the fact that $[i-2,j+2]$ is a $B$-badly approximable range,
  it follows that $z \ll 2^{(j-i)/2}$. Since the denominators of
  consecutive convergents of~$\beta$ are always coprime, we also have
  $\gcd(q_{r},q_r + q_{r-1}) = 1$, which concludes the proof
  of~(\ref{bohr_inclusion}).

  The bound~(\ref{cardinality_bohr}) follows immediately,
  and~(\ref{average_bohr}) is obtained by positivity from
  Lemma~\ref{av_phi}, using the bounds on~$x$, $y$ and~$z$:
  \[
    \sum_{N\in\mathcal{N}(i,j)}\frac{N}{\varphi(N)} \leq \sum_{n\in
      P(x,y,z)}\frac{n}{\varphi(n)} \ll z^2+z\log(zy)\ll
    2^{j-i}+j2^{(j-i)/2}.
  \]
\end{proof}

We are now in the position to prove property~(\ref{eqn05}). 

\begin{lemma}\label{QIA}
{Let as before \[
  A_{p,i} = \Bigl[p\beta - \frac{c}{2^{i+1}}, p\beta +
    \frac{c}{2^{i+1}}\Bigr],\quad
  E_i= \bigcup_{2^i \leq p < 2^{i+1}}A_{p,i}.
\]
    and let $U_k,V_k$ as before.}
  For any $k$ sufficiently large, we have
  \[
    \sum_{i,j \in [U_k,V_k]} \lambda(E_i \cap E_j) \ll \Bigl(\sum_{i
      \in [U_k,V_k]} \lambda(E_i)\Bigr)^2,
  \]
  where the implied constant depends only on $B$ and $\delta$.
\end{lemma}

\begin{proof}
  It suffices to prove that
  \[
    \sum_{i,j \in [U_k,V_k]} \lambda(E_i \cap E_j) \ll 1,
  \]
  where the implied constant depends only on $B$ and $\delta$, in view
  of~(\ref{accum_meas}).
  
  Observe for ${2^i \leq p < 2^{i+1} \leq 2^j \leq p' < 2^{j+1}}$ the straightforward estimate
  \[
    \lambda(A_{p,i} \cap A_{p',j}) \ll {\frac{1}{2^j}}
    \mathds{1}_{\left[\lVert (p-p')\beta \rVert \leq {\frac{1}{2^i}}\right]}.
  \]
  
  Therefore, we obtain 
  \[
    \sum_{U_k \leq i \leq j \leq V_k} \lambda(E_i\cap E_j) \ll
    \sum_{U_k \leq j \leq V_k} \frac{1}{{2^j}} \sum_{i \leq j} C(i,j)
  \]
  where
  \[
    C(i,j) = \Bigl|\Bigl\{ (p,p')\,\mid\, p \leq {2^{i+1}},\ p' \leq
    {2^{j+1}}\,\mid\, \lVert (p-p')\beta \rVert \leq \frac{1}{{2^{i}}}
    \Bigr\}\Bigr|.
  \]
  {With $\mathcal{N}(i,j)$ as in \eqref{bohr_def},} we thus have

  \[
    C(i,j) \leq \sum_{N \in \mathcal{N}(i,j{+1})} |\{p \leq {2^{i+1}}\,\mid\,
    p + N\text{ is prime}\}|.
  \]

  Note that by the choice of $U_k \leq i \leq j \leq V_k \leq 2U_k$, the parameter~$N$ here
  satisfies
  \[N\leq {2^{j+1} \leq 2^{V_k}\leq (2^{U_k})^2 \leq (2^i)^2}.\] Thus Proposition~\ref{pr-primes}, (2) implies that
  \[
    C(i,j) \leq \sum_{N \in \mathcal{N}(i,j+1)} |\{p \leq {2^{i+1}}\,\mid\,
    p + N\text{ is prime}\}| \ll \sum_{N \in
      \mathcal{N}(i,j+1)}\frac{N}{\varphi(N)}\frac{2^{i}}{i^2} \ll
    \frac{2^i}{V_k^2} \sum_{N \in
      \mathcal{N}(i,j+1)}\frac{N}{\varphi(N)},
  \]
  where we used the inequalities $U_k \leq i \leq j \leq V_k \ll U_k$.
  
  We now split the range of $(U_k \leq i \leq j \leq V_k)$ into two
  parts: when $j -i < 2 \log j$, we use the classical pointwise
  upper-bound $\frac{N}{\varphi(N)} \ll \log \log N \ll \log {V_k}$ and
  \eqref{cardinality_bohr} to deduce
    \[
      \sum_{N \in \mathcal{N}(i,j+1)}\frac{N}{\varphi(N)} \ll 2^{j-i}
      \log {V_k}.
    \]
    
  For $j -i \geq 2 \log j$ we apply Lemma \ref{av_phi}
  to deduce 
  \[
    \sum_{N \in \mathcal{N}(i,j+1)}\frac{N}{\varphi(N)} \ll 2^{j-i}.
  \]
  
  Thus, it follows that
  \begin{align*}
    \sum_{U_k \leq i \leq j \leq V_k} \lambda(E_i\cap E_j) 
    &\ll \sum_{U_k \leq j \leq V_k} \frac{1}{2^{j}} \sum_{i \leq j}
      C(i,j) 
    \\&\ll \frac{1}{V_k^2}\sum_{U_k \leq j \leq V_k} \sum_{i \leq j} 1 +
    \frac{\log V_k}{V_k^2}\sum_{U_k \leq j \leq V_k}
    \sum_{j - 2\log V_k \leq i \leq j} 1
    \\&\ll 1 + \frac{V_k (\log V_k)^2}{V_k^2} \ll 1,
  \end{align*}
  which finishes the proof.
\end{proof}

We have thus checked all assumptions of Lemma~\ref{bhv_lem}, and this
concludes the proof of {Theorem~\ref{prime_twisted}}.

\subsection{Proof of Theorem \ref{th-deterministic}}\label{ssec-prime-twisted}

We now prove {Theorem~\ref{th-deterministic}}.  By {Theorem
\ref{prime_twisted}}, for any $\beta$ which is i.b.a. and any $c > 0$, we know that
\[
  \lambda\Bigl( \Bigl\{x \in [0,1)\,\mid\, {\Big\vert} x -
    {\{p\beta\}} {\Big\vert} < \frac{c}{p} \text{
    for infinitely many } p \Bigr\}\Bigr)=1.
\]

Using the trivial estimate
\[
  {\Bigl\lvert} \frac{\fl{p^2\beta}\pmod p}{p} - \{p\beta\} {\Bigr\rvert}
  \leq \frac{1}{p},
\]
 an application of
Cassels's Lemma (as in Corollary \ref{cor-cassels-2}) proves the statement.
        
\subsection{Counterexample}\label{sec-counterexample}

We prove here the final part of Theorem~\ref{th-deterministic}.  As
mentioned in the introduction, there exist i.b.a numbers of all diophantine exponents, so that no
condition on the diophantine exponent suffices to automatically rule
out that $\lfloor p^2\beta \rfloor \pmod p$ or
$\lfloor pn_p\beta\rfloor \pmod p$ satisfies
$\lambda(\mathcal{A}(\bm{a})) = 1$.

However, we show now that some condition is needed, since if $\beta$
is {``always very Liouville''} (meaning that all
partial quotients grow at least with a certain rate), then the property fails.

More precisely, let $\beta$ be any (irrational) number whose continued
fraction has the property that the convergent denominators $(q_k)_{k\geq 1}$
satisfy the condition
\[
  q_{k+1} \geq \exp(q_k)
\]
for $k\geq 1$. Recall that~$\bm{b}(\beta)$ is defined by
$b_p=\lfloor p^2 \beta \rfloor \pmod p$.  We will show that
$\lambda(\mathcal{A}(\bm{b},c))=0$ for some $c > 0$, hence in fact for
all~$c>0$ by Cassels's Lemma (Corollary~\ref{cor-cassels}).

To do this, we claim that
\begin{equation}\label{eq-counter}
  \lambda\Bigl(\bigcup_{\substack{q_k\leq p <
      q_{k+1}}}\left[\{p\beta\} - \tfrac{c}{p}, \{p\beta\} +
    \tfrac{c}{p}\right]\Bigr) \ll \frac{\log \log q_k}{\log q_k}
\end{equation}
for all~$k\geq 1$.  Since under the growth assumption on $q_k$, the
sequence $(q_k)_{k}$ grows (at least) doubly exponentially, we have
\[
  \sum_{k \geq 1} \frac{\log \log q_k}{\log q_k} < \infty,
\]
and thus the fact that $\lambda(\mathcal{A}(\bm{b},c))=0$ follows from
an application of the first Borel--Cantelli Lemma.

To prove~(\ref{eq-counter}), we will split 
the range of $p$ between $q_k$ and~$q_{k+1}$ into three subranges.

\par
\textbf{(Small primes)}: by the Mertens estimates, we have
\[
  \lambda\Bigl(\bigcup_{\substack{q_k\leq p \leq q_k\log q_k}}
  \left[\{p\beta\} - \tfrac{c}{p}, \{p\beta\} +
    \tfrac{c}{p}\right]\Bigr) {\ll} \sum_{q_k \leq p \leq q_k \log q_k}
  \frac{1}{p} \ll \frac{\log \log q_k}{\log q_k}.
\]
\par
\textbf{(Large primes)}:  by the Mertens estimates again, we have
\[
  \lambda\Bigl(\bigcup_{\substack{q_{k+1}/q_k^2\leq p \leq q_{k+1}}}
  \left[\{p\beta\} - \tfrac{c}{p}, \{p\beta\} +
    \tfrac{c}{p}\right]\Bigr) {\ll} \sum_{q_{k+1}/q_k^2 \leq p \leq
    q_{k+1}} \frac{1}{p} \ll \frac{\log q_k}{\log q_{k+1}} \ll
  \frac{1}{\log q_k}.
\]
\par
\textbf{(Medium primes)}: for the remaining primes, we use the
elementary estimate
\begin{equation}\label{da_fact}
  \beta = \frac{p_k}{q_k} + (-1)^k \frac{1}{q_kq_{k+1}} +
  O\left(\frac{1}{a_{k+1}q_kq_{k+1}}\right).
\end{equation}

Fix a residue class $j \pmod{q_k}$. Then note that
\begin{align*}
  \bigcup_{\substack{q_k\log q_k \leq p
  \leq \frac{q_{k+1}}{q_k^2}\\ p \equiv j \pmod{q_k}}}
  \left[\{p\beta\} - \tfrac{c}{p}, \{p\beta\} + \tfrac{c}{p}\right]
  &\subseteq \bigcup_{\substack{q_k\log q_k \leq n \leq
    \frac{q_{k+1}}{q_k^2}\\ n \equiv j \pmod{q_k}}}
  \left[\{n\beta\} - \tfrac{c}{n}, \{n\beta\} + \tfrac{c}{n}\right]
  \\
  &\subseteq \left[\frac{jp_k}{q_k} - \frac{2}{q_k^3} -
    \frac{c}{q_k\log q_k}, \frac{jp_k}{q_k} + \frac{2}{q_k^3}
    + \frac{c}{q_k\log q_k} \right],
\end{align*}
where we used \eqref{da_fact} in the last step. It follows that
\begin{align*}
  \lambda\Bigl( \bigcup_{\substack{q_k\log q_k \leq p \leq
  \frac{q_{k+1}}{q_k^2}}} \left[\{p\beta\} - \tfrac{c}{p},
  \{p\beta\} + \tfrac{c}{p}\right] \Bigr)
  & \leq \sum_{j =
    0}^{q_k-1}\lambda\Bigl(\bigcup_{\substack{q_k\log q_k \leq p \leq
    \frac{q_{k+1}}{q_k^2}\\ p \equiv j \pmod {q_k}}}
  \left[\{p\beta\} - \tfrac{c}{p}, \{p\beta\} +
  \tfrac{c}{p}\right]\Bigr) \\
  &\ll \frac{1}{\log q_k}.
\end{align*}

The combination of these three ranges proves~(\ref{eq-counter}), and
finishes the proof.

\section{Hausdorff dimensions of exceptional sets}
\label{sec-HD}

We actually prove the following stronger statement.

\begin{theorem}\label{Thm_HD}
  Let $\beta$ be badly approximable and
  $\bm{x} = (x_n)_{n \in \mathbb{N}}$ with
  $x_n = \{n\beta\} + O(1/n)$. Then there exists a continuous function
  $f = f_{\beta}$ with $\lim_{\delta \to 0}f_{\beta}(\delta) = 1$ such
  that there exists $N_0 = N_0(\beta,\delta,\bm{x})$ with
  \[\dim_H\left\{\alpha \in [0,1)\,\mid\,
      \forall n \geq N_0,\ \lvert \alpha - x_n \rvert >
      \frac{\delta}{n}\right\} \geq f_{\beta}(\delta).\]
\end{theorem}

{We recall that for a prime number~$p$, we
denote by $n_p\geq 1$ the positive integer such that $p$ is the
$n_p$-th prime in increasing order. Further, we recall the sequences $\bm{a}(\beta) = (a_p)_{p}$ and~$\bm{b}(\beta) = (b_p)_{p}$
  defined by
  \begin{equation}\label{recall_pnnp}
    a_p = \lfloor pn_p\beta\rfloor \pmod p,\quad\quad b_p = \lfloor
    p^2\beta \rfloor \pmod p.
  \end{equation}
Note that 
\[\left(\frac{\lfloor pn_p\beta\rfloor \pmod {p}}{p}\right)_{p \in \mathbb{P}} =\left(\frac{\lfloor np_n\beta\rfloor \pmod {p_n}}{p_n}\right)_{n \in \mathbb{N}} = \{n\beta\} +
O(1/p_n) = \{n\beta\} +
O(1/n).\]
Further, since $\frac{\lfloor p^2\beta\rfloor \pmod {p}}{p} = \{p\beta\} + O(1/p)$ is simply a subsequence of $x_n = \{n\beta\} + O(1/n)$, we obtain immediately the following
Corollaries:
}

\begin{corollary}
  {Let $\bm{b}(\beta)$ be as in \eqref{recall_pnnp} with $\beta$ badly approximable.} Then we have
  \[
    \lim_{c \to 0}\dim_H([0,1] \setminus \mathcal{A}(\bm{b},c)
    ) = 1.
  \]
\end{corollary}

\begin{corollary}
  {Let $\beta$ be badly approximable.} Then we have
  \[
    \lim_{c \to 0} \dim_H\left\{\alpha \in [0,1)\,\mid\, \left\lvert
        \alpha - \frac{ \lfloor np_n\beta\rfloor \pmod {p_n}}{p_n}
      \right\rvert < \frac{c}{n} \text{ only finitely often} \right\}
    = 1.
  \]

  In particular, for any $c > 0$ and $\bm{a}(\beta)$ as in \eqref{recall_pnnp}, we have
  \[
    \dim_H([0,1] \setminus \mathcal{A}(\bm{a},c) ) = 1.
  \]
\end{corollary}

\begin{remark}
  {These} results can be related to those known from the theory of twisted Diophantine
  approximation: Bugeaud, Harrap, Kristensen, and Velani~\cite{BHKV10}
  proved that the set
  \[
    \left\{x \in [0,1]\,\mid\, \liminf_{n \to \infty}n\lVert n\beta
      - x \rVert > 0 \right\}
  \]
  has full Hausdorff dimension for \textit{every} irrational $\beta$,
  and Tseng~\cite{T09} proved the stronger statement that the above set
  is a $1/8$-winning set in the sense of Schmidt games \cite{S66}. 
  In contrast to the Lebesgue-measure setting (since there is no analogue to Corollary \ref{cor-cassels-2}),
  these results do not imply the result we obtain here. Nevertheless, it
  might well be that the methods used in their articles are flexible enough to extend to
  our situation.
  However, we opted for a short independent proof which works for
  every badly approximable $\beta$.
\end{remark}

We will make use of the generalized Cantor set as developed by
Badziahin and Velani \cite[Theorem 4]{BV11}. For a detailed
description see \cite{BV11}; here we only note that the classical
middle third Cantor set falls in the framework with $R_N = 3$ and
$r_{N,N} = 1, r_{N,N-k} = 0$ for all $k < N$. We will only use the
special case where $\bm{R} = (R,R,\ldots)$ for the proof here, the
more general case will be necessary for the greedy algorithm.

\begin{lemma}\label{lem_rich_cantor}
  Let $\bm{R} = (R_N)_{N \in \mathbb{N}}$ and
  $\bm{r} = (r_{m,n})_{m \leq n} \in \mathbb{N}$ with $R_N \geq 4$ and
  $I \subseteq [0,1]$ closed. Denote the generalized Cantor set by
  \[K(I,\bm{R},r) = \bigcap_{N = 1}^{\infty}\bigcup\limits_{J \in
      \mathcal{J}_N} J.\] Here $\mathcal{J}_N$ is the collection of
  intervals at stage $N$ where each $\mathcal{J}_{N-1}$ was partitioned
  into $R_N$ subintervals of equal length, with
  $\mathcal{J}_0 = \{I\}$. For each
  $J \in \mathcal{J}_{n-k-1}, 0 \leq k \leq n-1$, we delete at most
  $r_{n,k}$ many intervals $I_n\subset J$ at stage $n$. Given
  $K(I,\bm{R},r)$, suppose that $R_N \geq 4$ for all $N \in \N$ and
  that for all $N \in \N$,
  \begin{equation}\label{deletion_cond}
    \sum_{k = 0}^n \left(r_{N-k,n} \prod_{i = 1}^k
      \left(\frac{4}{R_{N-1}}\right)\right)\leq \frac{R_N}{4}.
  \end{equation}
  
  Then we have
  \[
    \dim_H(K(I,\bm{R},r)) \geq \liminf_{N \to \infty} (1 - \log_{R_N}
    (2)).
  \]
\end{lemma}

\begin{proof}[Proof of Theorem \ref{Thm_HD}]
  For the proof of Theorem \ref{Thm_HD}, we build generalized Cantor
  sets in order to apply Lemma \ref{lem_rich_cantor}. For this, let $K$
  be such that $\lvert x_n - \{n\beta\}\rvert < \frac{K}{n}$ and
  $c = c(\beta) > 0$ be such that $n\lVert n\beta\rVert \geq c$ for
  all $n \in \N$ (such a $c$ exists by the assumption of $\beta$ being
  badly approximable).  Since the sets under consideration are nested
  in $\delta$, it suffices to look at $\delta = \frac{1}{R^2}$ where
  $R \in \mathbb{N}$ is sufficiently large in terms of $K$ and $c$.
  We will construct a generalized Cantor set
  $K([0,\tfrac{1}{R}],\bm{R},\bm{r})$ that contains the exceptional
  set with
  \begin{equation}\label{r_n_bound}
    r_{N,N-k} \ll_{K,c} \begin{cases} R &\text{ if } k = 1,\\
      0 &\text{ if } k \neq 1
    \end{cases}\end{equation}

where $\alpha \in K([0,\tfrac{1}{R}],\bm{R},\bm{r})$ implies $\lvert \alpha - x_n| > \frac{\delta}{n}$ for $n \geq N_0(R)$.
  Since $r_{N,N-k} = 0$ if $k \neq 1$, \eqref{deletion_cond}
  reduces to the question of whether $r_{N,N-1} \leq
  \frac{R^2}{16}$. Thus the statement follows after the construction from
  Lemma \ref{lem_rich_cantor} when $R$ is chosen sufficiently large by
  choosing $f_{\beta}(\delta) = 1 - \log_{1/\delta}(2)$ and
  $\delta \to 0$.

  We now let $I = [0,\frac{1}{R}]$. Given a stage $N \geq 2$, we will
  ensure that for all $J \in \mathcal{J}_N$, we have
  \begin{equation}\label{cond_stage_N}
    \lvert \alpha - x_n\rvert > \frac{\delta}{n}, \quad \alpha \in J,
    \quad R^{N-1} \leq n < R^N.
  \end{equation}
  
  Assuming \eqref{cond_stage_N} to hold for each $N \geq 2$, then for
  each
  $x \in \bigcap\limits_{N = 1}^{\infty}\bigcup\limits_{J \in
    \mathcal{J}_N} J$ we have no solution to
  $\lvert \alpha - x_n\rvert < \frac{\delta}{n}$ with $n > R$, which implies that
  $\bigcap\limits_{N = 1}^{\infty}\bigcup\limits_{J \in \mathcal{J}_N}
  J$ is contained in the exceptional set.

  To achieve \eqref{r_n_bound} and \eqref{cond_stage_N}
  simultaneously, we do the following: For any
  $J \in \mathcal{J}_{N-1}$ fixed, let $\mathcal{I}_{N+1}$ be the
  collection of the subdivision of $J$ into $R^2$ intervals $I_{N+1}$
  with length $\frac{1}{R^{N+2}}$. We now upper-bound the number of
  $I_{N+1} \subset J$ such that there exists $R^{N} \leq n < R^{N+1}$ with
  $S_n = [x_n - \frac{\delta}{n},x_n + \frac{\delta}{n}]$ satisfying
  $S_n \cap I_{N+1} \neq \emptyset$. In order to do so, we observe
  that
  \[
    |\{R^{N} \leq n < R^{N+1}\,\mid\, x_n \in J\}| \leq
    \Bigl|\left\{R^{N} \leq n < R^{N+1}\,\mid\, \{n\beta\} \in J +
      \left[-\tfrac{K}{R^N}, \tfrac{K}{R^N}\right]\right\}\Bigr|.
  \]
  Since $\lvert J\rvert = \frac{1}{R^N}$, the set
  $J + \left[-\tfrac{K}{R^N}, \tfrac{K}{R^N}\right]$ is an interval of
  length $\leq \frac{2K+1}{R^N}$. Note that since $\beta$ is badly
  approximable, we have
  \[
    \min_{R^{N} \leq n\neq m < R^{N+1}}\lVert n\beta - m\beta \rVert
    \geq \frac{c}{R^{N+1}},
  \]
  so we obtain
  \[
    |\{R^{N} \leq n < R^{N+1}\,\mid\, x_n \in J\}| \leq
    \frac{2K+1}{R^{N}}\frac{R^{N+1}}{c} + 1 \leq R \frac{2K+1}{c} + 1.
  \]
  Furthermore note that for each fixed $n$, we have that $S_n$ is an
  interval of length at most $\frac{\delta}{R^N} = \frac{1}{R^{N+2}}$
  and thus can intersect at most three intervals $I_{N+1}$. Therefore,
  we obtain
  \[
    \Bigl| \Bigl\{ I_{N+1} \subseteq J\,\mid\, \Bigl(\bigcup_{R^N \leq n
        < R^{N+1}} S_n\Bigr) \cap I_{N+1} \neq \emptyset\Bigr\} \Bigr|
    \ll_{K,c} R.
  \]
  
  Removing these intervals and keeping all remaining ones in step $N+1$
  proves that we can choose $r_{N,N-1} \ll_{K,c} R$, which proves the
  statement.
\end{proof}

Finally, using similar again the generalized Cantor set idea, we prove here the final
statement of Theorem~\ref{thm_greedy}.

\begin{theorem}\label{thm_greedy_HD}
  Let $\uple{g}= (a_p)_{p \in \mathbb{P}}$ be the sequence constructed
  in Theorem~\ref{thm_greedy}.  Then for $c$ sufficiently small, we
  have
  \[\dim_H\left([0,1) \setminus \mathcal{A}(\bm{g},c) \right) = 1.\]
\end{theorem}

\begin{proof}
	Our construction will only consider the interval $I = [1/4,3/4]$ and we will again construct 
    a generalized Cantor set $K(I,\bm{R},r)$. In contrast to the proof of Theorem \ref{Thm_HD}, we will have $r_{N,k} = 0$ for all $k < N$, but 
    $(R_n)_{n \in \N}$ will not be a constant sequence, but determined by the number of steps taken in each iteration of the Greedy algorithm.\\
	
	Let $p_j^{-}$ the biggest prime $p$ such that the Greedy algorithm is in iteration $j$ and \[I_{p} := \left(\tfrac{a_p}{p} - \tfrac{c}{p},\tfrac{a_p}{p} + \tfrac{c}{p}\right) \subseteq [0,1/4].\] Analogously, let
	$p_j^{+}$ the smallest prime in iteration $j$ such that $I_{p} \subseteq [3/4,1]$. Next, we define define $R_n$ recursively in the following way:
    
	Let $C_{n-1} := \prod_{j \leq n-1}R_j$ be given (since $R_1,\ldots,R_{n-1}$ have been constructed). We then claim that for $n$ sufficiently large and $C_{n-1} < p_n^{-}$
	there exists $k_n \in \mathbb{N}$ such that:
	
	\begin{itemize}
		\item[(i)] $p_{n}^{+} = o(k_nC_{n-1}), \quad k_nC_{n-1} = o(p_{n+1}^{-}),$
		\item[(ii)] $\liminf_{n \to \infty} k_n = \infty$. 
	\end{itemize}
Indeed, by the construction of the algorithm and Mertens estimate 
(since $\sum_{X < p < Y} \frac{1}{p} \approx 1/2$ implies $Y \approx X^{\sqrt e}$ for $X,Y$ sufficiently large, so in particular $X = o(Y)$),
we obtain $\lim_{n \to \infty}\frac{p_n^{+}}{p_n^{-}} = \infty, \lim_{n \to \infty}\frac{p_n^{-}}{p_{n+1}^{-}} = \infty$.
Setting $k_n := \left\lfloor \frac{\sqrt{p_{n+1}^{-}p_{n}^{+}}}{C_{n-1}}\right\rfloor$, we observe that
\[k_nC_{n-1} = \sqrt{p_{n+1}^{-}p_{n}^{+}} + O(C_{n-1}) = \sqrt{p_{n+1}^{-}p_{n}^{+}} + O(p_n^{-}) = \sqrt{p_{n+1}^{-}p_{n}^{+}} + o(p_n^{+}).\] 
This proves (i) and since $\frac{\sqrt{p_{n+1}^{-}p_{n}^{+}}}{C_{n-1}} \geq \frac{p_{n}^{+}}{p_n^{-}}\to \infty$, (ii) follows as well.
	So now we set $R_n = k_n$ provided $n$ is sufficiently large (for formal reasons, we set $R_n =1$ for the first boundedly many $n$ to ensure $C_{n-1} < p_n^{-}$).\\
	
	Having the parameters $(R_n)$ established, we fix $n$ sufficiently large and assume an interval ${J}_{n-1} \in \mathcal{J}_{n-1}$ to be constructed. Note that by construction 
	$|J_{n-1}| = \frac{1}{2C_{n-1}}$. Next, we partition $J_{n-1}$ into $k_n$ many subintervals $J'$ of length $\frac{1}{2k_nC_{n-1}}$
	We will now cross out all subintervals $J'$ that intersect with an interval $I_p = \left(\tfrac{a_p}{p}- \tfrac{c}{p},\tfrac{a_p}{p} + \tfrac{c}{p}\right)$
	that is generated in the $n$-th iteration of the greedy algorithm. Note that this is equivalent to considering $I_p$ with 
	$p_{n}^{-} \leq p \leq p_n^{+}$. We see that such an interval $I_p$ can only cover $\leq \frac{2c}{p}k_nC_{n-1} + 2$ many
	subintervals $J'$ of $J$. Note that $p < p_n^{+} = o(k_nC_{n-1})$, so for large enough $n$, $\frac{c}{p}k_nC_{n-1} > 2$, so we cover at most
	$\frac{3c}{p}k_nC_{n-1}$ many subintervals of $J$. Now observe that most primes do not intersect $J$ at all: Let $p_{J}^{-},p_{J}^{+}$ be the
	smallest respectively largest prime $p$ (in iteration $n$ of the Greedy algorithm) such that $I_p \cap J = \emptyset$. Then the construction of the algorithm implies that
	
	\[\frac{\sum_{p_{J}^{-} \leq p \leq p_{J}^{+}} \frac{1}{p}}{|J|} \in \left[\tfrac{1}{2},5\right],\]
	provided $n$ is sufficiently large: This follows from $a_p \geq a_q + 1/p, p_n^- \leq q<p \leq p_n^{+}$ (see proof of Lemma \ref{lem_greedy}) and the observation that there can only be $O(1)$ many $p \in [p_n^-,p_n^{+}]$ where $I_p$ intersects both $J$ and $I \setminus J$, so the contribution of these $p$ is negligible since $\frac{\tfrac{1}{p}}{|J|} \leq \frac{2k_{n-1}C_{n-2}}{p_n^{-}} = o(1)$ by (i).	
	Thus the number of removed subintervals $J' \subseteq J$ is bounded from above by
	\[\sum_{p_{J}^{-} \leq p \leq p_{J}^{+}} \frac{3c}{p}k_nC_{n-1} < 30c|J|R_nC_{n-1} = 15cR_n.\]
	Therefore, if $c < 1/60$ and $n$ sufficiently large, we remove at most $R_n/4$ subintervals at stage $n$, while in the remaining $\geq 3R_n/4$ many intervals we avoid all intervals $I_p$ where $p$ is considered in the $n$-th iteration of the algorithm. In other words, this ensures that
	\begin{equation}
		r_{n,k} = \begin{cases}
			R_n/4 \text{ if } k = n,\\
			0 \text{ otherwise.}
		\end{cases}
	\end{equation}
	 This proves \eqref{deletion_cond}. 
	 By construction, if $x \in K(I,\bm{R},r)$, then the inequality $\left\lvert \frac{a_p}{p} - x\right\rvert < \tfrac{c}{p}$ can only have solutions for $p$ in the first (boundedly many) iterations of the algorithm, so in particular, the inequality has only finitely many solutions.  
	 Since $\liminf_{n \to \infty} R_n = \liminf_{n \to \infty} k_n = \infty$, an application of Lemma \ref{lem_rich_cantor} proves full Hausdorff dimension, provided $c < 1/60$.
\end{proof}

\bibliographystyle{plain}
\bibliography{bibliography.bib}

\begin{thebibliography}{10}

\bibitem{ABH23}
C.~Aistleitner, B.~Borda, and M.~Hauke.
\newblock On the metric theory of approximations by reduced fractions: a
  quantitative {K}oukoulopoulos-{M}aynard theorem.
\newblock {\em Compos. Math.}, 159(2):207--231, 2023.

\bibitem{BV11}
D.~Badziahin and S.~Velani.
\newblock Multiplicatively badly approximable numbers and generalised {C}antor
  sets.
\newblock {\em Adv. Math.}, 228:2766--2796, 2011.

\bibitem{BHV24}
V.~Beresnevich, M.~Hauke, and S.~Velani.
\newblock Borel--{C}antelli, zero-one laws and inhomogeneous
  {D}uffin--{S}chaeffer.
\newblock {\em Preprint. \url{arXiv:2406.19198}}, 2024.

\bibitem{BV08}
V.~Beresnevich and S.~Velani.
\newblock A note on zero-one laws in metrical {D}iophantine approximation.
\newblock {\em Acta Arith.}, 133:363--374, 2008.

\bibitem{BHKV10}
Y.~Bugeaud, S.~Harrap, S.~Kristensen, and S.~Velani.
\newblock On shrinking targets for {$\mathbb{Z}^m$} actions on tori.
\newblock {\em Mathematika}, 56:193--202, 2010.

\bibitem{C50}
J.W.S. Cassels.
\newblock Some metrical theorems in {D}iophantine approximation. {I}.
\newblock {\em Math. Proc. Cambridge Philos. Soc.}, 46:209--218, 1950.

\bibitem{Cas57}
J.W.S. Cassels.
\newblock {\em An introduction to {D}iophantine approximation}, volume No. 45
  of {\em Cambridge Tracts in Mathematics and Mathematical Physics}.
\newblock Hafner Publishing Co., New York, 1972.
\newblock Facsimile reprint of the 1957 edition.

\bibitem{chow}
S.~Chow.
\newblock Bohr sets and multiplicative {D}iophantine approximation.
\newblock {\em Duke Math. J.}, 167(9):1623--1642, 2018.

\bibitem{dv86}
H.G. Diamond and J.D. Vaaler.
\newblock Estimates for partial sums of continued fraction partial quotients.
\newblock {\em Pacific J. Math.}, 122(1):73--82, 1986.

\bibitem{FK16}
M.~Fuchs and D.~H. Kim.
\newblock On {K}urzweil's 0-1 law in inhomogeneous {D}iophantine approximation.
\newblock {\em Acta Arith.}, 173:41--57, 2016.

\bibitem{Harman_1998}
G.~Harman.
\newblock {\em Metric number theory}.
\newblock Oxford Univ. Press, 1998.

\bibitem{ant}
H.~Iwaniec and E.~Kowalski.
\newblock {\em Analytic number theory}, volume~53 of {\em Amer. Math. Soc.
  Colloq. Publ.}
\newblock American Mathematical Society, Providence, RI, 2004.

\bibitem{kowalski_ergodic}
E.~Kowalski.
\newblock Unmotivated ergodic averages.
\newblock {\em preprint, \url{arXiv:2309.13576}}, 2023.

\bibitem{K55}
J.~Kurzweil.
\newblock On the metric theory of inhomogeneous diophantine approximations.
\newblock {\em Studia Math.}, 15:84--112, 1955.

\bibitem{lubotzky-meiri}
A.~Lubotzky and C.~Meiri.
\newblock Sieve methods in group theory {I}: {P}owers in linear groups.
\newblock {\em J. Amer. Math. Soc.}, 25:1119--1148, 2012.

\bibitem{R20}
F.A. Ram\'irez.
\newblock Khintchine's theorem with random fractions.
\newblock {\em Mathematika}, 66:178--199, 2020.

\bibitem{S66}
W.M. Schmidt.
\newblock On badly approximable numbers and certain games.
\newblock {\em Trans. Amer. Math. Soc.}, 123:178--199, 1966.

\bibitem{tao_blog}
T.~Tao.
\newblock Continued fractions, {B}ohr sets, and the {L}ittlewood conjecture.
\newblock
  \url{https://terrytao.wordpress.com/2012/01/03/continued-fractions-bohr-sets-and-the-littlewood-conjecture/}.

\bibitem{T09}
J.~Tseng.
\newblock Badly approximable affine forms and {S}chmidt games.
\newblock {\em J. Number Theory}, 129:3020--3025, 2009.

\end{thebibliography}

\end{document}